\documentclass[12pt]{amsart}
\usepackage{cases}
\usepackage{txfonts}
\usepackage{latexsym}
\usepackage{amsthm}
\usepackage{mathrsfs}
\usepackage{amssymb, amsmath}

\def\opn#1#2{\def#1{\operatorname{#2}}} % to make operators
\opn\chara{char} \opn\length{\ell}
\opn\projdim{proj\,dim} \opn\injdim{inj\,dim} \opn\rank{rank}
\opn\depth{depth} \opn\grade{grade} \opn\height{height}
\opn\embdim{emb\,dim} \opn\codim{codim}

\opn\Tr{Tr} \opn\bigrank{big\,rank}
\opn\superheight{superheight}\opn\lcm{lcm}
\opn\trdeg{tr\,deg}%
\opn\reg{reg} \opn\lreg{lreg}
\opn\div{div} \opn\Div{Div} \opn\cl{cl} \opn\Cl{Cl}
\opn\Spec{Spec} \opn\Supp{Supp} \opn\supp{supp} \opn\Sing{Sing}
\opn\Ass{Ass}
\opn\Ann{Ann} \opn\Rad{Rad} \opn\Soc{Soc}
\opn\Ker{Ker} \opn\Coker{Coker} \opn\Im{Im} \opn\Hom{Hom}
\opn\Tor{Tor} \opn\Ext{Ext} \opn\End{End} \opn\Aut{Aut} \opn\id{id}

\opn\nat{nat}
\opn\pff{pf}%   \pf exists already
\opn\Pf{Pf} \opn\GL{GL} \opn\SL{SL} \opn\mod{mod} \opn\ord{ord}
\opn\aff{aff} \opn\con{conv} \opn\relint{relint} \opn\st{st}
\opn\lk{lk} \opn\cn{cn} \opn\core{core} \opn\vol{vol}
\opn\gr{gr}

\def\pot#1#2{#1[\kern-0.28ex[#2]\kern-0.28ex]}

\opn\dirlim{\underrightarrow{\lim}}
\opn\invlim{\underleftarrow{\lim}}
\def\qed{\ifhmode\textqed\fi
   \ifmmode\ifinner\quad\qedsymbol\else\dispqed\fi\fi}
\def\textqed{\unskip\nobreak\penalty50
    \hskip2em\hbox{}\nobreak\hfil\qedsymbol
    \parfillskip=0pt \finalhyphendemerits=0}
\def\dispqed{\rlap{\qquad\qedsymbol}}

\opn\ini{in} \opn\inm{inm} \opn\Sym{Sym}

\begin{document}
\title{On large gaps between zeros of the Riemann zeta-function}
\author{ShaoJi Feng and XiaoSheng Wu}
\date{}
\address {Academy of Mathematics and Systems Science, Chinese Academy of Sciences,
Beijing 100190, P. R.
China.}
\email {fsj@amss.ac.cn}
\address {Academy of Mathematics and Systems Science, Chinese Academy of Sciences,
Beijing 100190, P. R.
China.}
\email {xswu@amss.ac.cn}

\begin{abstract} Assuming the Generalized Riemann Hypothesis(GRH), we show
that infinitely often consecutive non-trivial zeros of the Riemann
zeta-function differ by at least 3.072 times the average spacing.

\end{abstract}
\maketitle

\section{Introduction}
Let $\zeta(s)$ denote the Riemann zeta-function. We denote the
non-trivial zeros of $\zeta(s)$ as $\rho=\beta+i\gamma$. Let
$\gamma\le\gamma'$ denote consecutive ordinates of the zeros of
$\zeta(s)$. The von Mangoldt formulate (see \cite{Tichmarsh}) gives
\begin{align}
   N(T)=\frac T{2\pi}\log\frac T{2\pi e}+O(\log T),\notag
\end{align}
where $N(T)$ is the number of zeros of $\zeta(s), s=\sigma+it$ in
the rectangle $0\le\sigma\le1, 0\le t\le T$. Hence, the average size
of $\gamma'-\gamma$ is $2\pi/\log\gamma$. In 1973, by studying the
pair correlation of the zeros of the Riemann zeta-function,
Montgomery \cite{Montgomery1} suggested that there exists arbitrarily large and
small gaps between consecutive zeros of $\zeta(s)$. That is to say
\begin{align}
   \lambda=\lim\sup(\gamma'-\gamma)\frac{\log\gamma}{2\pi}=\infty \ \ \ \
   \mathrm{and} \ \ \ \
   \mu=\lim\inf(\gamma'-\gamma)\frac{\log\gamma}{2\pi}=0\notag,
\end{align}
where $\gamma$ runs over all the ordinates of the zeros of the
$\zeta(s)$.

In this article, we focus on the large gaps and assume the Generalized Riemann Hypothesis (GRH) is true.
This conjecture states that the non-trivial zeros of the Dirichlet $L$-functions are on the Re(s)=1/2 line. We obtain\\\\
{\bf Theorem 1.1}\ \ \textit{If the Generalized Riemann Hypothesis is
true, then $\lambda>3.072$}.\\

Unconditionally, selberg remarked in \cite{Selberg} that he could prove
$\lambda>1$. Assuming RH, Mueller \cite{Mueller} showed that $\lambda>1.9$, and
later, by a different approach, Montgomery and Odlyzko \cite{Montgomery2} obtained
$\lambda>1.9799$. This result was then improved by Conrey, Ghosh,
and Gonek \cite{Conrey1} who obtained $\lambda>2.337$ assuming RH and
$\lambda>2.68$ in \cite{Conrey2} assuming GRH. Recently, by making use of the
Wirtinger inequality, Hall \cite{Hall} proved that there exist infinity many
large gaps between the zeros on the critical line of the Riemann
zeta-function greater than 2.63 times the average spacing of the
zeros of Riemann zeta-function. This result implies $\lambda>2.63$
on RH. Assuming Riemann Hypothesis, H.M.Bui, M.B.Milinovich and N.Ng proved $\lambda>2.69$ in \cite{Bui2} and we obtained $\lambda>2.7327$ in \cite{Feng2}. On GRH, N.Ng \cite{Ng} proved in 2006 that $\lambda>3$
and this result was improved to $\lambda>3.033$ by Bui \cite{Bui1} in 2009.

The works of \cite{Bui1}, \cite{Conrey2}, \cite{Ng} are based on the following idea of
J.Mueller \cite{Mueller}. Let $H:\mathbb{C}\rightarrow \mathbb{R}_{\geq0}$ be continuous and
consider the associated functions
\begin{align}
\label{1.1}
   \mathcal {M}_1(H,T)=\int_1^TH(\frac12+it)dt,\tag{1.1}
\end{align}
\begin{align}
\label{1.2}
   m(H,T;\alpha)=\sum_{T<\gamma<2T}H(\frac12+i(\gamma+\alpha))dt,\tag{1.2}
\end{align}
\begin{align}
\label{1.3}
   \mathcal {M}_2(H,T;c)=\int_{-c/L}^{c/L}m(H,T;\alpha)d\alpha,\tag{1.3}
\end{align}
where $L=\log{\frac T{2\pi}}$. Note that
\begin{align}
\label{1.4}
   \frac{\mathcal {M}_2(H,2T;c)-\mathcal {M}_2(H,T;c)}{\mathcal
   {M}_1(H,2T)-\mathcal {M}_1(H,T)}<1 \tag{1.4}
\end{align}
implies $\lambda>\frac c\pi$.

Mueller applied this idea with $H(s)=|\zeta(s)|^2$. Let $A(s)$
denote a Dirichlet polynomial
\begin{align}
\label{1.5}
   A(s)=\sum_{n\le y}a(n)n^{-s}.\tag{1.5}
\end{align}
On RH, Conrey et al. used $H(s)=|A(s)|^2$ with
$a(n)=d_{2.2}(n),y=T^{1-\epsilon}$ and obtained $\lambda>2.337$.
Here $d_r(n)$ is the coefficient of $n^{-s}$ in the Dirichlet series
$\zeta(s)^r$. Later, assuming GRH, they applied (\ref{1.4}) to
$H(s)=|\zeta(s)A(s)|^2$ with $a(n)=1$ and
$y=(T/2\pi)^{1/2-\epsilon}$ and obtained $\lambda>2.68$. By
considering a more general coefficients $a(n)$, N.Ng \cite{Ng} proved
$\lambda>3$. Actually, N.Ng chose $H(s)=|\zeta(s)A(s)|^2$ with
$A(s)$ has coefficients $a_r(n)=d_r(n)p(\frac{\log n}{\log y})$. In \cite{Bui1}, H. M. Bui
chose $H(s)=A_1(s)+\zeta(s)A_2(s)$, where $A_1(s),A_2(s)$ are Dirichlet series defined by (\ref{1.5})
with coefficients $a_{r_1}=d_{r_1}(n)p_1(\frac{\log n}{\log y})$ and $a_{r_2}=d_{r_2}(n)p_2(\frac{\log n}{\log y})$, and proved $\lambda>3.033$.

In this article, we choose $H(s)=|\zeta(s)A(s)|^2$,
 where $A(s)$ is defined by (\ref{1.5}) with $y=T^{\frac12-\epsilon}$ and the coefficients
\begin{align}
\label{1.6}
   a(n)=d_r(n)P_0\bigg(\frac{\log n}{\log y}\bigg)+d_r^*(n)P_2\bigg(\frac{\log n}{\log
   y}\bigg),\tag{1.6}
\end{align}
for $P_0, P_2$ are polynomials and $r\in \mathbb{N}$. Here,
\begin{align}
\label{1.7}
   d_r^*(n)=\frac1{\log^2y}\Lambda*\Lambda*d_r(n),\tag{1.7}
\end{align}
where $*$ is the convolution and $\Lambda$ is the Mongoldt function.

Let \begin{align}
   \chi(s)=2^s\pi^{s-1}\sin\frac12s\pi\Gamma(1-s)\notag
\end{align}
and
\begin{align}
   Z(t)=(\chi(\frac12+it))^{-\frac12}\zeta(\frac12+it).\notag
\end{align}
It's well known that $Z'(t)$ has a zero between the consecutive
zeros of Riemann zeta-function. Since
\begin{align}
   |\chi(s)|=1 \ \ \ \ \mathrm{and} \ \ \ \ \frac{\chi'(\frac12+it)}{\chi(\frac12+it)}\sim\log\frac{t}{2\pi},\notag
\end{align}
we have
\begin{align}
   |Z'(t)|\sim\bigg|\zeta(\frac12+it)\bigg|\bigg|-\frac12\log\frac t{2\pi}+\frac{\zeta'}{\zeta}(\frac12+it)\bigg|.\notag
\end{align}
From (\ref{1.7}), it's easy to see
\begin{align}
   \zeta^r(s)\big(\frac{\zeta'}{\zeta}(s)\big)^2=(\log y)^2\sum_{n=1}^\infty\frac{d_r^*(n)}{n^s}\notag
\end{align}
for $s=\sigma+it$ with $\sigma>1$.
Hence, our choice of $H(s)$ may be seen as a kind of approximation to
$Z(t)^{2r-2}Z'(t)^4(\log y)^{-4}$.

We now come to the precise result. We define several functions that
will appear in the following. Given $P_0, P_2$ are polynomials and
$u\in\mathbb{Z}_{\ge0}$, we define
\begin{align}
\label{1.8}
   Q_{i,u}(x)=\int_0^1\theta^uP_i(x+\theta(1-x))d\theta,\tag{1.8}
\end{align}
for $i=0,2$. Given $\eta\in\mathbb{R}$ and
$\overrightarrow{n}=(n_1,n_2,n_3,n_4)\in(\mathbb{Z}_{\ge0})^4$, we
define
\begin{align}
\label{1.9}
   l_{i_1,i_2}(\overrightarrow{n})=\int_0^1x^{r^2+n_1+n_2-1}(1-x)^{2r+n_3+n_4}Q_{i_1,r+n_3-1}(x)Q_{i_2,r+n_4-1}(x)dx\tag{1.9}
\end{align}
\begin{align}
\label{1.10}
   k_{i_1,i_2}(\overrightarrow{n})=&\int_0^1\int_0^{1-x}
   x^{r+n_1-1}(\eta^{-1}-x)^{n_2}y^{r^2+n_3-1}(1-y)^{r+n_4}\notag\\
   &\cdot P_{i_1}(x+y)Q_{i_2,r+n_4-1}(y)dydx.\tag{1.10}
\end{align}
For
$\overrightarrow{n}=(n_1,n_2,n_3,n_4,n_5)\in(\mathbb{Z}_{\ge0})^5$,
we define
\begin{align}
\label{1.11}
   h_{i_1,i_2}(\overrightarrow{n})=&-\eta^{-1}\text{B}(n_5+1,r+n_4-1)l_{i_1,i_2}(n_1,n_2,n_3-1,n_4+n_5)\notag\\
   &+\text{B}(n_5+1,r+n_4-1)l_{i_1,i_2}(n_1,n_2,n_3,n_4+n_5)\notag\\
   &+\text{B}(n_5+1,r+n_4)l_{i_1,i_2}(n_1,n_2,n_3-1,n_4+n_5+1),\tag{1.11}
\end{align}
where $\text{B}(m,n)$ is the Beta function. For $r\ge1$, we define the
constants
\begin{align}
   i_1=i_1'+i_1'',\ \ \ \ \ \ \ \ i_2=i_2'+i_2''\notag
\end{align}
\begin{align}
   a_r=&\prod_{p}((1-p^{-1})^{r^2}\sum_{m=0}^\infty(\frac{\Gamma(r+m)}{\Gamma(r)m!})^2p^{-m}),\notag
\end{align}
\begin{align}
\label{1.12}
   b_r(i_1',i_2')=\sum_{\tau=0}^{\min(i_1',i_2')}C_{i_1'}^\tau
   C_{i_2'}^\tau\tau! r^{i_1'+i_2'-2\tau}\tag{1.12}
\end{align}
\begin{align}
\label{1.13}
   c_r(i_1',i_2',i_1'',i_2'')=\frac{C_{i_1}^{i_1'}C_{i_2}^{i_2'}b_r(i_1',i_2')}{(r^2+i_1'+i_2'-1)!(r+i_1''-1)!(r+i_2''-1)!},\tag{1.13}
\end{align}
with $C_{m}^n$ is the
binomial coefficient. For $n\ge-2$, we also define
\begin{align}
\label{1.14}
   \Omega_{r}(i_2'',n)=\left\{ \begin{aligned}
           &(-1)^{n+1}C_r^{n+2} \ \ \ \ \ \ \ \ \ \ \ \ \ \ \ \ \ \ \ \ \ \ \ \ \ \ \ \ \ \ \ \ \ \ \ \ \ \ \ \ \ \ \ \ \ \text{for} \ \ \ \ i_2''=0\notag\\
           &\sum_{j'=-2}^{\min(r-2,n)}(-1)^{j'+1}C_r^{j'+2}\Delta(n-j')\ \ \ \ \ \ \ \ \ \ \ \ \ \ \ \ \ \ \ \ \text{for} \ \ \ \ i_2''=1\notag\\
           &\sum_{j'=-2}^{\min(r-2,n)}(-1)^{j'+1}C_r^{j'+2}\sum_{j_1+j_2=n-j'}\Delta(j_1)\Delta(j_2)\ \ \ \ \ \text{for} \ \ \ \
           i_2''=2\notag
   \end{aligned} \right.\tag{1.14}
\end{align}
with $\Delta$ given by
\begin{align}
\label{1.15}
\left\{ \begin{aligned}
         \Delta(0) &= 1\notag\\
                  \Delta(j)&=-1,\ \ \ \ \ \ \ \ \ \    \text{for} \ \ \ \ \ \ j\ge1.
                          \end{aligned} \right.\tag{1.15}
                          \end{align}
                          Since
\begin{align}
   \sum_{j=0}^{r}(-1)^{j+1}C_r^{j}P(j)=0\notag
\end{align}
for any polynomial $P(j)$ on $j$, it's not difficult to see
\begin{align}
   \Omega_{r}(i_2'',n)=0 \ \ \ \ \ \ \text{for} \ \ \ \ \ n>r-2.\notag
\end{align}
From this definitions, we can present our result for
$\mathcal{M}_1(H,T)$ and $m(H,T;\alpha)$.\\\\
{\bf Theorem 1.2}\ \ \textit{Let $y=(\frac T{2\pi})^\eta$ with $0<\eta<1/2$, we have
\begin{align}
\label{1.16}
   \mathcal{M}_1(H,T)\sim a_{r+1}T(\log y)^{(r+1)^2}
   \sum_{i_1'+i_1''=0,2}\sum_{i_2'+i_2''=0,2}
   c_r(i_1',i_2',i_1'',i_2'')\hat{l}(\eta,r,i_1',i_2',i_1'',i_2'')\tag{1.16}
\end{align}
as $T\rightarrow\infty$, where
\begin{align}
\label{1.17}
   \hat{l}(\eta,r,i_1',i_2',i_1'',i_2'')=&\eta^{-1}l_{i_1,i_2}(i_1',i_2',
   i_1'',i_2'')\notag\\
   &-l_{i_1,i_2}(i_1',i_2',i_1''+1,i_2'')\notag\\
   &-l_{i_1,i_2}(i_1',i_2',i_1'',i_2''+1).\tag{1.17}
\end{align}
This is valid up to an error term which is $O(L^{-1})$ smaller than
the main term}.\\\\
{\bf Theorem 1.3}\ \ \textit{Suppose $r\in \mathbb{N}$ and $y=(\frac T{2\pi})^\eta$ with $\eta<1/2$. The
Generalized Riemann Hypothesis implies
\begin{align}
\label{1.18}
   m(H,T;\alpha)\sim&\frac{a_{r+1}TL^{(r+1)^2+1}}{\pi}\mathrm{R}e\sum_{j=0}^\infty
   z^j\eta^{j+(r+1)^2+1}\bigg(\frac{\hat{h}(r,j,\eta)}{j!}+\hat{k}(r,j,\eta)\bigg)\notag\\
   &+\frac{L}{2\pi}\mathcal{M}_1(H,T),\tag{1.18}
\end{align}
where $z=i\alpha L, |z|\ll1$,
\begin{align}
\label{1.19}
   \hat{h}(r,j,\eta)=&\sum_{i_1'+i_1''=0,2}\sum_{i_2'+i_2''=0,2}c_r(i_1',i_2',i_1'',i_2'')\notag\\  &\cdot(rh_{i_1,i_2}(i_1',i_2',i_1''+1,i_2''+1,j)+i_2''(r+i_2''-1)h_{i_1,i_2}(i_1',i_2',i_1''+1,i_2'',j+1)),\tag{1.19}
\end{align}
\begin{align}
\label{1.20}
   \hat{k}(r,j,\eta)=&\sum_{i_1'+i_1''=0,2}\sum_{i_2'+i_2''=0,2}c_r(i_1',i_2',i_1'',i_2'')\notag\\
   &\cdot\sum_{n=-2}^{\min(r-2,j)}\frac{\Omega_{r}(i_2'',n)(r+i_2''-1)!}
   {(j-n)!(r+i_2''+n+1)!} k_{i_1,i_2}(i_1'',j-n,i_1'+i_2',i_2''+n+2).\tag{1.20}
\end{align}
This result is valid up to an error term
$O_{\epsilon,r}(TL^{(r+1)^2}+T^{1/2+\eta+\epsilon})$}.

From Theorem 1.2 and Theorem 1.3, an argument similar to N.Ng \cite{Ng}
deduce that
\begin{align}
   \frac{\mathcal{M}_2(H,2T;c)-\mathcal{M}_2(H,T;c)}{\mathcal{M}_1(H,2T)-\mathcal{M}_1(H,T)}=f_r(c)+O(\epsilon),\notag
\end{align}
where
\begin{align}
\label{1.21}
   f_r(c)=\frac1D\sum_{j=0}^\infty\frac{(-1)^jc^{2j+1}}{2^{2j}} \bigg(\frac{\hat{h}(r,2j,\frac12)}{(2j+1)!}+\frac{\hat{k}(r,2j,\frac12)}{2j+1}\bigg)+\frac
   c\pi+O(\epsilon)\tag{1.21}
\end{align}
and
\begin{align}
   D:=\pi\sum_{i_1'+i_1''=0,2}\sum_{i_2'+i_2''=0,2}
   c_r(i_1',i_2',i_1'',i_2'')\hat{l}_{i_1,i_2}(\eta,r,i_1',i_2',i_1'',i_2'').\notag
\end{align}
It's known that $f_r(c)<1$ implies $\lambda\ge\frac c\pi$. We may
compute (\ref{1.21}) for various choices of $r$ and $P_0(x), P_2(x)$.
Choosing $c=3.072\pi, r=2$ and $P_0(x)=x^{30}, P_2(x)=-31.4x^{165}$, we
compute the sum
\begin{align}
   D^{-1}\sum_{j=0}^{30}\frac{(-1)^jc^{2j+1}}{2^{2j}}\bigg(\frac{\hat{h}(r,2j,\frac12)}{(2j+1)!}+\frac{\hat{k}(r,2j, \frac12)}{2j+1}\bigg)+\frac
   c\pi=0.999846...\notag
\end{align}
by \textit{Mathematic}. On the other hand, we may bound the terms $j>30$. For $P_0(x)=x^{30}, P_2(x)=-31.4x^{165}$, it's easy to see $|Q_{i,u}(x)|\le32$ on [0,1] and $l_{i_1,i_2}(\overrightarrow{n})\le32^2$.
So, a direct calculation gives that $\widehat{h}(\overrightarrow{n})\le64\times4\times32^2$ and hence
\begin{align}
   \bigg|\frac1D\sum_{j>30}^\infty\frac{(-1)^jc^{2j+1}}{2^{2j}}\frac{\widehat{h}(r,2j,\frac12)}{(2j+1)!}\bigg| &\le\frac{262144c}{D}\sum_{j>30}\frac{(c/2)^{2j}}{(2j+1)!}\notag\\
   &\le\frac{262144c}{D}\sum_{j>30}e^{-2j(\log(2j)-(\log(c/2)-1))}\notag\\
   &<\frac{262144c}{D}\frac{e^{-60(\log(60)-\log(c/2)-1)}}{2\log(60)-\log(c/2)-1}<10^{-20},\notag
\end{align}
where we have applied $n!>(n/e)^n$. A similar calculation establishes that
\begin{align}
   \bigg|\frac1D\sum_{j>J}^\infty\frac{(-1)^jc^{2j+1}}{2^{2j}} \frac{\widehat{k}(r,2j,\frac12)}{(2j+1)}\bigg|<10^{-20}.\notag
\end{align}
Thus, we conclude that $f_2(3.05\pi)<1$ and
establish Theorem 1.1. If we let $r=2$ and $P_0(x)=x^{30},
P_2(x)=0$, we get $f_2(3\pi)=0.999481...$, which accords with the result of
N.Ng \cite{Ng}.

We have deduced Theorem 1.1 from Theorem 1.2 and Theorem 1.3. Hence,
the rest of the article will be devoted to establishing the result
of Theorem 1.2 and Theorem 1.3. From a similar argument to the part
4 of N.Ng \cite{Ng}, we note
\begin{align}
\label{1.22}
   m(H,T;\alpha)\sim2\mathrm{Re}I+\frac L{2\pi}\mathcal {M}_1(H,T)\tag{1.22}
\end{align}
with an error term $O(L^{-1})$ smaller. Here,
\begin{align}
\label{1.23}
   I=\sum_{k\le
   y}\frac{a(k)}{k}\sum_{j\le{kT}{2\pi}}b(j)e(-j/k)+O(yT^{\frac12+\epsilon}),\tag{1.23}
\end{align}
\begin{align}
\label{1.24}
   b(j)=-\sum_{{hmn=j}\atop{h\le
   y}}a(h)d(m)\Lambda(n)n^{i\alpha}.\tag{1.24}
\end{align}
Hence, to prove Theorem 1.2 and Theorem 1.3, it's sufficient to evaluate
$\mathcal {M}_1$ and $I$. We will evaluate $\mathcal {M}_1$ in
section 4 and $I$ in section 5.
\section{Some notation and definitions}
Throughout this article we shall employ the notation
\begin{align}
\label{2.1}
   [m]_y:=\frac{\log m}{\log y}\tag{2.1}
\end{align}
for $m,y>0$, and we appoint that $p,p_i,q,q_j$ always denote primes for
$i,j\ge1$. The sum
\begin{align}
   \sum_{a_1+\cdots+a_m\ge D}\ \ \ \ \text{and} \ \ \ \ \sum_{a_1+\cdots+a_m=D}\notag
\end{align}
are always over all entire arrays $(a_1,a_2,\cdots, a_m)$
with $a_i\ge0$. In addition, we define
$j_{\tau}(n)$ and $\sigma_r(n)$ as in N.Ng \cite{Ng},
\begin{align}
\label{2.2}
   j_{\tau}(n)=\prod_{p\mid n}(1+O(p^{-\tau}))\tag{2.2}
\end{align}
for $\tau>0$ and the constant in the $O$ is fixed and independent of
$\tau$ and 
\begin{align}
\label{2.3}
   \sigma_r(n)=\prod_{p^\lambda\|n}d_r(p^\lambda)H_{\lambda,r}(p^{-1})\tag{2.3}
\end{align}
with
\begin{align}
   H_{\lambda,r}(x):=\lambda
   x^{-\lambda}\int_0^xt^{\lambda-1}(1-t)^{r-1}dt.\notag
\end{align}
Here, $p^\lambda\|n$ means $p^\lambda|n$ and $p^{\lambda+1}\nmid n$.
A simple calculation by part integration shows that
\begin{align}
   H_{\lambda,1}(x)=1, \ \ \ \ \ H_{\lambda,2}(x)=1-\frac{\lambda}{\lambda+1}x,\notag
\end{align}
and for $r\ge3$,
\begin{align}
   &H_{\lambda,r}(x)\notag\\
   =&(1-x)^{r-1}+\sum_{i=1}^{r-2}\frac{(r-1)\cdots(r-i)}
   {(\lambda+1)\cdots(\lambda+i)}x^i(1-x)^{r-i-1}+\frac{(r-1)!\lambda!}{(\lambda+r-1)!}x^{r-1}.\notag
\end{align}
From this, it's easy to see $H_{\lambda,r}(x)$ is a polynomial of
$x$ with $H_{\lambda,r}(0)=1$, and all the coefficients of the
polynomial are $O(1)$. Here, the constant of the $O$ is only decided
by $r$. So, we have
\begin{align}
\label{2.4}
   \sigma_r(p_1\cdots
   p_i)=r^i+O(\sum_{\tau=1}^i\frac1{p_\tau}),\tag{2.4}
\end{align}
\begin{align}
\label{2.5}
   \sigma_r(p_1\cdots p_im)\ll\sigma_r(p_1\cdots
   p_i)\sigma_r(m)+O(\sigma_r(m)\sum_{\tau=1}^i\frac1{p_\tau})\tag{2.5}
\end{align}
with the constant of $O$ is only decided by $r$ and $i$, for $m,i\ge1$ are integers. We now also invoke several properties of
$d_r$ which we apply repeatedly as follow:
\begin{align}
\label{2.6}
   &\sum_{m\le x}d_r(m)m^{-1}\ll\log^rx,\notag\\
   &\sum_{m\le x}d_r(m)^2m^{-1}\ll\log^{r^2}x.\tag{2.6}
\end{align}
\section{Some lemmas}
 In this section, we present some lemmas that will be used in the following.\\\\
{\bf Lemma 3.1}\ \ (Mertens Theorem).
\textit{\begin{align}
\label{3.1}
   \sum_{p\le y}\frac {\log p}{p}=\log y+O(1).\tag{3.1}
\end{align}}\\
{\bf Lemma 3.2}\ \ \textit{For positive integers $m_1$,$m_2$ and $n$,
\begin{align}
\label{3.3}
   &\sum_{{p_1p_2\cdots p_{m_1}\mid n}\atop{q_1q_2\cdots q_{m_2}\mid n}}\mu^2(p_1\cdots p_{m_1})\log p_1\cdots\log p_{m_1}
   \mu^2(q_1\cdots q_{m_2})\log q_1\cdots\log q_{m_2}\nonumber\\
   =&\sum_{k=0}^{\min (m_1,m_2)}\mathcal {C}_{m_1}^k\mathcal{C}_{m_2}^kk!
   \sum_{p_1\cdots p_{m_1+m_2-k}\mid n}\mu^2(p_1\cdots p_{m_1+m_2-k})\notag\\
   &\cdot\log ^2p_1\cdots\log ^2p_k\log p_{k+1}\cdots\log p_{m_1+m_2-k}\tag{3.2}
\end{align}\\
where $p$ and $q$ runs over prime numbers, $ \mathcal {C}_{m}^{k}$
is the binomial coefficient}.

This Lemma is a generalization of Lemma 2.3 in Feng \cite{Feng1}.\\\\
{\bf Lemma 3.3}\ \ \textit{Let $a_i\ge1$ be integers for $1\le i\le m$,
$F(x)\ll M$ on $[1,y]$ be continuous,
\begin{align}
\label{3.4}
   &\sum_{p_1\cdots p_m\le
   y}\frac{\log^{a_1}p_1\cdots\log^{a_m}p_m}{p_1\cdots
   p_m}\int_1^{\frac y{p_1\cdots p_m}}\frac{F(p_1\cdots
   p_mx)\log^kx}{x}dx\notag\\
   =&\frac{k!\prod_{i=1}^m(a_i-1)!}{(\sum_{i=1}^ma_i+k)!}\int_1^y\frac{F(x)(\log
   x)^{\sum_{i=1}^ma_i+k}}{x}dx+O(M(\log
   y)^{\sum_{i=1}^ma_i+k})\tag{3.3}
\end{align}}
\textit{Proof}. By Lemma 3.1 and Abel summation, we may express the
right side of (\ref{3.4}) as the expression in Lemma 9 of Feng \cite{Feng1}. Then, an argument similar to
the proof of Lemma 9 Feng \cite{Feng1} establishes the Lemma.\\\\
{\bf Lemma 3.4}\ \ (Conrey \cite{Conrey3} Lemma 3). \textit{Suppose that
$A_j(s)=\sum_{n=1}^\infty\alpha_j(n)n^{-s}$ is absolutely convergent
for $\sigma>1$, for $1\le j\le J$, and that
\begin{align}
\label{3.5}
   A(s)=\sum_{n=1}^\infty\frac{\alpha(n)}{n^s}=\prod_{j=1}^JA_j(s).\tag{3.4}
\end{align}
Then for any positive integer $d$,
\begin{align}
\label{3.6}
   \sum_{n=1}^\infty\frac{\alpha(dn)}{n^s}=\sum_{d_1\cdots
   d_J=d}\prod_{j=1}^J\bigg(\sum_{(n,d_1\cdots
   d_{j-1})=1}^\infty\frac{\alpha_j(nd_j)}{n_s}\bigg).\tag{3.5}
\end{align}}
{\bf Lemma 3.5}\ \ (N.Ng \cite{Ng} Lemma 5.3). \textit{Let $(h,k)=1$ and $k=\prod
p^\lambda>0$. For $\alpha\in \mathbb{R}$ and $\sigma>1$ define}
\begin{align}
\label{3.7}
  \mathcal
  {Q}(s,\alpha,h/k)=-\sum_{m,n=1}^\infty\frac{d(m)\Lambda(n)}{m^sn^{s-i\alpha}}e\bigg(\frac{-mnh}{k}\bigg).\tag{3.6}
\end{align}
\textit{Then ${Q}(s,\alpha,h/k)$ has a meromorphic continuation to the
entire complex plane, If $\alpha\neq0$, ${Q}(s,\alpha,h/k)$ has}

(i) \textit{at most a double pole at $s=1$ with same principal part as}
\begin{align}
   k^{1-2s}\zeta^2(s)\bigg(\frac{\zeta'}{\zeta}(s-i\alpha)-\mathcal
   {G}(s,\alpha,k)\bigg),\notag
\end{align}
\textit{where}
\begin{align}
\label{3.8}
   \mathcal
   {G}(s,\alpha,k)=\sum_{p\mid k}\log
   p\bigg(\sum_{a=1}^{\lambda-1}p^{a(s-1+i\alpha)}+\frac{p^{\lambda(s-1+i\alpha)}}{1-p^{-s+i\alpha}}-\frac1{p^{s-i\alpha}-1}\bigg);\tag{3.7}
\end{align}

(ii) \textit{a simple pole at $s=1+i\alpha$ with residue}
\begin{align}
   -\frac1{k^{i\alpha}\phi(k)}\zeta^2(1+i\alpha)\mathcal
   {R}_k(1+i\alpha)\notag
\end{align}
\textit{where}
\begin{align}
\label{3.9}
   \mathcal
   {R}_k(s)=\prod_{p^\lambda\parallel
   k}(1-p^{-1}+\lambda(1-p^{-s}(1-p^{s-1})).\tag{3.8}
\end{align}
\textit{Moreover, on GRH,} $\mathcal {Q}(s,\alpha,h/k)$ \textit{is regular in
$\delta>1/2$ except for these two poles}.\\\\
{\bf Lemma 3.6}\ \ \textit{Assume GRH. Let $y=(\frac T{2\pi})^\eta$ where
$0<\eta<1/2$, $k\in \mathbb{N}$ with $k\le y$, and $\alpha\in \mathbb{R}$. Set}
\begin{align}
\label{3.10}
   \mathcal {Q}^*(s,\alpha,k)=\sum_{j=1}^\infty b(j)j^{-s}e(-j/k)\ \
   \ (\sigma>1),\tag{3.9}
\end{align}
\textit{where}
\begin{align}
   b(j)=-\sum_{{hmn=j}\atop{h\le
   y}}(d_r(h)P_1([h]_y)+d_r^*(h)P_2([h]_y))d(m)\Lambda(n)n^{i\alpha}.\notag
\end{align}
\textit{Then} $\mathcal {Q}^*(s,\alpha,k)$ \textit{has an analytic continuation to
$\sigma>1/2$ except possible poles at $s=1$ and $1+i\alpha$.
Furthermore,}
\begin{align}
   \mathcal {Q}^*(s,\alpha,k)=O(y^{\frac12}T^\epsilon)\notag
\end{align}
\textit{where $s=\sigma+it, 1/2+L^{-1}\le\sigma\le1+L^{-1}, |t|\le T,
|s-1|>0.1$, and $|s-1-i\alpha|>0.1$}.\\\\
\textit{Proof}. From the definition of $b(j)$, we may denote $\mathcal {Q}^*(s,\alpha,k)=\mathcal
{Q}_1^*(s,\alpha,k)+\mathcal {Q}_2^*(s,\alpha,k)$ with obvious
meaning and prove both parts satisfy the Lemma. The proof of
$\mathcal {Q}_1^*(s,\alpha,k)$ is given by Lemma 5.6 of N.Ng \cite{Ng}. We
can prove $\mathcal {Q}_2^*(s,\alpha,k)$ similarly to Lemma 5.6 of
N.Ng \cite{Ng}. The only difference is we replace (5.9) of N.Ng \cite{Ng} with
\begin{align}
   B(s,d,z)=&\sum_{f_1f_2f_3f_4f_5f_6=d}A_1(s,f_1;z)A_2(s,f_2,f_1)A_2(s,f_3,f_1f_2)\notag\\
   &\times A_3(s,f_4,f_1f_2f_3)A_4(s,f_5,f_1f_2f_3f_4)A_4(s,f_6,f_1f_2f_3f_4f_5)\notag
\end{align}
by Lemma 3.4, where
\begin{align}
   A_1(s,f;z)&=\chi(f)\sum_{h\le
   y/f}\frac{\chi(h)d_r(fh)(fh)^z}{h^s},\notag\\
   A_2(s,f,r)&=\chi(f)L(s,\chi)\prod_{p\mid r}(1-\chi(p)p^{-s}),\notag\\
   A_3(s,f,r)&=-\sum_{(n,r)=1}\chi(fn)\Lambda(fn)(fn)^{i\alpha}n^{-s},\notag\\
   A_4(s,f,r)&=-\sum_{(n,r)=1}\chi(fn)\Lambda(fn)n^{-s}.\notag
\end{align}
It's obvious that $A_4(s,f,r)=A_3(s,f,r)$ for $\alpha=0$, so, the
other part of the proof is the same to Lemma 5.6 of N.Ng \cite{Ng}.\\\\
{\bf Lemma 3.7}. \textit{For $\alpha\in \mathbb{R}$ and $j\in \mathbb{Z}_{\ge0}$, we have}
\begin{align}
\label{3.11}
  \mathcal{G}^{(j)}(1,\alpha,k)=\sum_{p\mid k}p^{i\alpha}(\log
  p)^{j+1}+O(C_j(k))\tag{3.10}
\end{align}
\textit{where} $\mathcal{G}(s,\alpha,k)$ \textit{is defined by (\ref{3.8}) and
\begin{align}
   C_j(k)=\sum_{p\mid k}\frac{\log^j
   p}{p}+\sum_{p^i\|k,i\ge2}\alpha\log^jp.\notag
\end{align}
Moreover, for $x\le y$, we have
\begin{align}
\label{3.12}
   \sum_{h,k\le
   x}\frac{a(h)a(k)(h,k)}{hk}C_j(\frac{k}{(h,k)})\ll(\log
   x)^{r^2+2r},\tag{3.11}
\end{align}}
\textit{Proof}. We remark that (\ref{3.11}) is proven in Conrey \cite{Conrey2}. Recalling
the definition of $a(n)$, we may denote the left side of (\ref{3.12}) as
\begin{align}
   \sum_{h,k\le
   x}\frac{C_j(\frac{k}{(h,k)})(h,k)}{hk}\big(d_r(h)d_r(k)+d_r^*(h)d_r(k)+d_r(h)d_r^*(k)+d_r^*(h)d_r^*(k)\big).\notag
\end{align}
Thus, we express the left side of (\ref{3.12}) into four parts. The first
part accords with (\ref{3.12}) given by Lemma 5.7 N.Ng \cite{Ng} and we now prove
it's also available to the other three parts. We only give the proof
of the fourth part, since the other parts can be proven similarly.
The part we are considering is
\begin{align}
\label{3.13}
   &\sum_{h,k\le
   x}\frac{d_r^*(h)d_r^*(k)(h,k)}{hk}C_j(\frac{k}{(h,k)})\notag\\
   &\le\sum_{h,k\le
   x}\frac{d_r^*(h)d_r^*(k)}{hk}(C_j(k)+1)\sum_{a\mid(h,k)}\phi(a)\notag\\
   &\le\sum_{a\le x}\frac1{a}\sum_{h,k\le\frac
   xa}\frac{d_r^*(ah)d_r^*(ak)(C_j(ak)+1)}{hk},\tag{3.12}
\end{align}
where $\phi(n)$ is the number of numbers less than $n$ and prime to $n$.
Recalling that
\begin{align}
   d_r^*(n)=\frac1{\log^2y}\sum_{i_1,i_2=1}^\infty\sum_{p_1^{i_1}p_2^{i_2}\mid n}\log p_1\log p_2d_r\bigg(\frac
   n{p_1^{i_1}p_2^{i_2}}\bigg)\notag
\end{align}
and
\begin{align}
   C_j(ak)=\sum_{p\mid ak}\frac{\log^j
   p}{p}+\sum_{p^i\|ak,i\ge2}\alpha\log^jp,\notag
\end{align}
we find the sum in (\ref{3.13}) is
\begin{align}
\label{3.14}
   &\ll\frac1{\log^4x}\sum_{a\le x}\frac1a\sum_{h\le\frac
   xa}\frac{1}{h}\sum_{i_1,i_2,j_1,j_2=1}^\infty\sum_{p_1^{i_1}p_2^{i_2}\mid ah}\log p_1\log p_2d_r\bigg(\frac{ah}{p_1^{i_1}p_2^{i_2}}\bigg)\notag\\
   &\cdot\sum_{k\le\frac xa}\frac1k\sum_{q_1^{j_1}q_2^{j_2}\mid
   ak}\log q_1\log
   q_2d_r\bigg(\frac
   {ak}{q_1^{j_1}q_2^{j_2}}\bigg)\bigg(\sum_{p\mid ak}\frac{\log^j
   p}{p}+\sum_{p^i\|ak,i\ge2}\alpha\log^jp\bigg).\tag{3.13}
\end{align}
We divide the sum in (\ref{3.14}) into five parts by the number of
different elements in $\{p_1,p_2,q_1,q_2,p\}$. Not shortage of
general nature, we only prove the part with any two elements are
different here, for the other parts can be proven similarly. we find
the part consisted by the terms with any two elements in $\{p_1,p_2,q_1,q_2,p\}$ are different
in the sum of (\ref{3.14}) is
\begin{align}
   \ll&\frac1{\log^4x}\sum_{p\le x}\frac{\log^jp}{p^2}\sum_{i_1,i_2,j_1,j_2=1}^\infty\sum_{p_1^{i_1}\le
   x}\frac{\log p_1}{p_1^{i_1}}\sum_{p_2^{i_2}\le x}\frac{\log p_2}{p_2^{i_2}}\sum_{q_1^{j_1}\le x}\frac{\log q_1}{q_1^{j_1}}
   \sum_{q_2^{j_2}\le x}\frac{\log q_2}{q_2^{j_2}}\notag\\
   &\cdot\sum_{a\le x}\frac{d_r^2(a)}{a}\sum_{h\le
   x}\frac{d_r(h)}{h}\sum_{k\le x}\frac{d_r(k)}{k}\ll((\log
   x)^{r^2+2r}),\notag
\end{align}
for familiar formula
\begin{align}
   \sum_{i\ge2}\sum_{p^i\le x}\frac{\log^jp}{p^i}=O(1)\notag
\end{align}
with $\forall j\ge0$. Putting together the results establishes the
lemma.\\\\
{\bf Lemma 3.8}\ \ \textit{Suppose $r,n\in \mathbb{N}, 1\le x, n\le\frac T{2\pi}$,
and $F\in C^1([0,1])$. There exists an absolute constant
$\tau_0=\tau_0(r)$ such that}
\begin{align}
\label{3.15}
   \sum_{h\le x}\frac{d_r(nh)}{h}F([h]_x)=\frac{\sigma_r(n)(\log x)^r}{(r-1)!}
   \int_0^1\theta^{r-1}F(\theta)d\theta+O(d_r(n)j_{\tau_0}(n)L^{r-1}),\tag{3.14}
\end{align}
\textit{with $j_{\tau_0}(n)$ defined by (\ref{2.2}). Furthermore, suppose $a_i\ge1$ are integers for $1\le i\le m$, we have}
\begin{align}
\label{3.16}
   &\sum_{h\le x}\frac{F([h]_x)}{h}\sum_{p_{1}\cdots
   p_{m}\mid
   h}\log^{a_1}p_{1}\cdots\log^{a_m}p_{m}d_r\bigg(\frac{nh}{p_1\cdots
   p_{m}}\bigg)\notag\\
   =&\frac{\sigma_r(n)\prod_{i=1}^m(a_i-1)!(\log x)^{r+\sum_{i=1}^ma_i}}{(r+\sum_{i=1}^ma_i-1)!}
   \int_0^1\theta^{r+\sum_{i=1}^ma_i-1}F(\theta)d\theta\notag\\
   &+O(d_r(n)j_{\tau_0}(n)L^{r+\sum_{i=1}^ma_i-1})\tag{3.15}.
\end{align}
\textit{Proof}. The first identity (\ref{3.15}) is given by lemma 5.8 N.Ng \cite{Ng}.
Changing summation order and making the variable change
$h\rightarrow hp_1\cdots p_m$ yields the left side of (\ref{3.16})
\begin{align}
   =&\sum_{p_1\cdots
   p_m\le x}\frac{\log^{a_1}p_1\cdots\log^{a_m}p_m}{p_1\cdots
   p_m}\sum_{h\le x/p_1\cdots p_m}\frac{F([p_1\cdots
   p_mh]_x)}{h}d_r(nh)\notag\\
   =&\frac{\sigma_r(n)}{(r-1)!}\sum_{p_1\cdots
   p_m\le x}\frac{\log^{a_1}p_1\cdots\log^{a_m}p_m}{p_1\cdots
   p_m}\int_1^x(\log t)^{r-1}F([p_1\cdots
   p_mt]_x)\frac{dt}{t}\notag\\
   &+O\bigg(\sum_{p_1\cdots
   p_m\le x}\frac{\log^{a_1}p_1\cdots\log^{a_m}p_m}{p_1\cdots
   p_m}d_r(n)j_{\tau_0}(n)L^{r-1}\bigg).\notag
\end{align}
Then (\ref{3.16}) follows by Lemma 3.1 and Lemma 3.3.\\\\
{\bf Lemma 3.9}\ \ \textit{For $r,i\in \mathbb{N}$, $1\le i_1\le i$ and $g\in
C^1([0,1])$, we have}
\begin{align}
   &\sum_{n\le y}\frac{\phi(n)\sigma_r(p_1\cdots p_{i_1}n)\sigma_r(p_{i_1+1}\cdots
   p_in)}{n^2}g([n]_y)\notag\\
   =&\frac{r^ia_{r+1}(\log y)^{r^2}}{(r^2-1)!}\int_0^1\theta^{r^2-1}g(\theta)d\theta+O((\log y)^{r^2}\big(\sum_{\tau=1}^ip_{\tau}^{-1}+(\log
   y)^{-1}\big).\notag
\end{align}
\textit{Moreover, suppose $a_i\ge1$ are integers for $1\le i\le m$, let $1\le i_1<i_2<\cdots<i_{m_1}\le m$ and $1\le i'_1<i'_2<\cdots<i'_{m_2}\le m$ for $0\le m_1,m_2\le m$,
then}
\begin{align}
   &\sum_{n\le y}\frac{\phi(n)}{n^2}g([n]_y)\sum_{p_{1}\cdots
   p_{m}\mid
   h}\log^{a_1}p_{1}\cdots\log^{a_m}p_{m}\sigma_r\bigg(\frac{n}{p_{i_1}\cdots p_{i_{m_1}}}\bigg)\sigma_r\bigg(\frac{n}{p_{i'_1}\cdots p_{i_{m_2}'}}\bigg)\notag\\
   &\sim\frac{r^{2m-m_1-m_2}a_{r+1}\prod_{i=1}^m(a_i-1)!(\log y)^{r^2+\sum_{i=1}^ma_i}}{(r^2+\sum_{i=1}^ma_i-1)!}
   \int_0^1\theta^{r^2+\sum_{i=1}^ma_i-1}g(\theta)d\theta\notag
\end{align}
\textit{plus an error $O((\log y)^{-1})$ smaller.} \\\\
\textit{Proof}. We remark that the first identity is a generalization
of Lemma 5.9 (i) in N.Ng, and it can be proven similarly equal to
\begin{align}
   &\frac{\sigma_r(p_1\cdots p_{i_1})\sigma_r(p_{i_1+1}\cdots p_i)a_{r+1}(\log y)^{r^2}}{(r^2-1)!}
   \int_0^1\theta^{r^2-1}g(\theta)d\theta\notag\\
   &+O((\log y)^{r^2}(\sum_{\tau=1}^ip_{\tau}^{-1}+(\log
   y)^{-1}),\notag
\end{align}
then the first identity follows by (\ref{2.4}). The second identity can be
proven by the first identity with an argument as the proof of (\ref{3.16}) in Lemma 3.8.\\

We define $f(k)=\mathcal{R}_k(1+i\alpha)/\phi(k)$ and
$\mathcal{T}_{k;N}(\alpha)=\sum_{j=0}^N\mathcal{R}_k^{(j)}(1)(i\alpha)^j/j!$ with $\mathcal{R}_k(s)$ given by (\ref{3.9}).\\\\
{\bf Lemma 3.10} (N.Ng \cite{Ng} Lemma 5.11).\ \ \textit{For $l=\log x$,
$|\alpha|\ll(\log x)^{-1}$, $1\le x,m\le y$, $n$ square free and
$n\mid m$,  we have}
\begin{align}
\label{3.17}
   \sum_{k\le
   x}d_r(mk)f(nk)\ll\frac{d_r(m)j_{\tau_0}(m)l^r}{n^{1-\epsilon}},\tag{3.16}
\end{align}
\textit{where $\tau_0=1/3$ is valid}.\\\\
{\bf Lemma 3.11}\ \ \textit{Let $l=\log x, |\alpha|\ll(\log x)^{-1}$,
$\tau_0=1/3$ and $g\in C^1([0,1])$. We have}
\begin{align}
\label{3.18}
   &\sum_{k\le
   x}d_r(mk)g([k]_x)\frac{\mathcal{R}_{nk}^{(j)}(1)}{\phi(nk)}\notag\\
   =&\frac{\sigma_r(m)(-1)^jC_r^j(\log x)^{r+j}}{n(r+j-1)!}\int_0^1\theta^{r+j-1}g(\theta)\frac{dt}t
   +O\bigg(\frac{d_r(m)j_{\tau_0}(m)l^{r+j-1}}{n^{1-\epsilon}}\bigg)\tag{3.18}
\end{align}
\textit{and}
\begin{align}
\label{3.19}
   \sum_{k\le
   x}d_r(mk)\bigg(f(nk)-\frac{\mathcal{T}_{nk;r}(\alpha)}{\phi(nk)}\bigg)
   \ll|\alpha|^{r+1}l^{2r}\frac{d_r(m)j_{\tau_0}(v)}{n^{1-\epsilon}}.\tag{3.19}
\end{align}
\textit{Still, suppose $a_i\ge1$ are integers for $1\le i\le \tau$, $k=p_1\cdots
p_\tau k'$, then}
\begin{align}
\label{3.20}
   &\sum_{k\le
   x}\frac{g([k]_x)}{k}\sum_{p_1\cdots p_\tau\mid k}\log^{a_1}p_1\cdots\log^{a_\tau}p_{\tau}d_r(mk')\frac{nk'\mathcal{R}_{nk'}^{(j)}(1)}{\phi(nk')}\notag\\
   =&\frac{\sigma_r(m)(-1)^jC_r^j\prod_{i=1}^\tau(a_i-1)!(\log x)^{r+\sum_{i=1}^\tau a_i+j}}{(r+\sum_{i=1}^\tau a_i+j-1)!}
   \int_0^1\theta^{r+\sum_{i=1}^\tau a_i+j-1}g(\theta)d\theta\notag\\
   &+O\bigg(\frac{d_r(m)j_{\tau_0}(m)l^{r+\sum_{i=1}^\tau a_i+j-1}}{n^{1-\epsilon}}\bigg)\tag{3.20}
\end{align}

The identities (\ref{3.18}) and (\ref{3.19}) are given by N.Ng \cite{Ng}, and the identity (\ref{3.20})
can be proven by (\ref{3.18}) with an argument as the proof of (\ref{3.16}) in Lemma 3.8. \\\\
{\bf Lemma 3.12}\ \ \textit{Let $A(s)=\sum_{n\le y}\frac{a(n)}{n^s}$, where
$y=(\frac T{2\pi})^\eta$ and $\eta\in(0,\frac12)$. Then for $1\le
t\le T$,}
\begin{align}
\label{3.21}
   \int_0^t|\zeta A(\frac12+iu)|^2du=t\sum_{h,k\le
   y}\frac{a(h)a(k)(h,k)}{hk}\log\frac{t(h,k)^2e^{2\gamma-1}}{2\pi
   hk}+O(T),\tag{3.21}
\end{align}
\textit{here $\gamma$ is Euler's constant}.

This lemma is a special case of a formula of Balasubramanian, Conrey
and Heath-Brown \cite{Balas}.

\section{Evaluation of $\mathcal {M}_1$}
From (1.1) we recall that
\begin{align}
   \mathcal {M}_1(H,T)=\int_1^T|\zeta A(\frac12+it)|^2dt.\notag
\end{align}
Then by Lemma 3.12,
\begin{align}
   \mathcal {M}_1(H,T)=T\sum_{h,k\le
   y}\frac{a(h)a(k)(h,k)}{hk}\log\frac{T(h,k)^2e^{2\gamma-1}}{2\pi
   hk}+O(T).\notag
\end{align}
To estimate the sum we apply the M\"{o}bius inversion formula
\begin{align}
   f((h,k))=\sum_{{m\mid h}\atop{m\mid k}}\sum_{n\mid
   m}\mu(n)f(\frac mn),\notag
\end{align}
and obtain
\begin{align}
   \mathcal {M}_1(H,T)=T\sum_{h,k\le
   y}\frac{a(h)a(k)}{hk}\sum_{{m\mid h}\atop{m\mid k}}\sum_{n\mid m}\frac{\mu(n)m}{n}\log\frac{Te^{2\gamma-1}m^2}{2\pi
   n^2hk}+O(T).\notag
\end{align}
Changing the order of summation and replacing $h$ by $hm$, $k$ by
$km$, we find that
\begin{align}
   \mathcal {M}_1(H,T)=T\sum_{m\le y}\frac1m\sum_{n\mid m}\frac{\mu(n)}{n}\sum_{h,k\le
   y/m}\frac{a(mh)a(mk)}{hk}\log\frac{Te^{2\gamma-1}}{2\pi
   n^2hk}+O(T).\notag
\end{align}
We next replace the logarithm term by $\log(T/(2\pi hk))$ with an
error $O(\log n)$. A calculation shows that this $O(\log n)$ term
contributes $O(TL^{r^2+2r})$ in $\mathcal {M}_1(H,T)$. Since
$\sum_{n\mid m}\mu(n)n^{-1}=\phi(m)m^{-1}$ we deduce that
\begin{align}
   \mathcal {M}_1(H,T)=T\sum_{m\le y}\frac{\phi(m)}{m^2}\sum_{h,k\le
   y/m}\frac{a(mh)a(mk)}{hk}\log\frac{T}{2\pi
   hk}+O(TL^{r^2+2r}).\notag
\end{align}
Recalling the definition of $a(n)$, we denote
\begin{align}
   \mathcal {M}_1(H,T)=&T\sum_{m\le y}\frac{\phi(m)}{m^2}\sum_{h,k\le
   y/m}\frac{\log\frac{T}{2\pi
   hk}}{hk}\big(d_r(mh)d_r(mk)P_1([mh]_y)P_1([mk]_y)\notag\\
   +&d_r^*(mh)d_r(mk)P_2([mh]_y)P_1([mk]_y)\notag\\
   +&d_r(mh)d_r^*(mk)P_1([mh]_y)P_2([mk]_y)\notag\\
   +&d_r^*(mh)d_r^*(mk)P_2([mh]_y)P_2([mk]_y)\big)+O(TL^{r^2+2r})\notag\\
   =&\mathcal {M}_{11}+\mathcal {M}_{12}+\mathcal {M}_{13}+\mathcal
   {M}_{14}+O(TL^{r^2+2r})\notag
\end{align}
with obvious meaning. We now come to calculate $\mathcal {M}_{14}$.
Recalling the definition of $d_r^*(n)$ by (\ref{1.7}), we observe that
\begin{align}
   d_r^*(mh)=&\frac1{\log^2y}\sum_{i_1,i_2=1}^\infty\sum_{p_1^{i_1}p_2^{i_2}\mid mh}\log
   p_1\log p_2d_r\bigg(\frac {mh}{p_1^{i_1}p_2^{i_2}}\bigg),\notag
\end{align}
\begin{align}
   d_r^*(mk)=&\frac1{\log^2y}\sum_{j_1,j_2=1}^\infty\sum_{q_1^{j_1}q_2^{j_2}\mid mk}\log
   q_1\log q_2d_r\bigg(\frac {mk}{q_1^{j_1}q_2^{j_2}}\bigg).\notag
\end{align}
We may replace $d_r^*(mh)$ and $d_r^*(mk)$ in $\mathcal{M}_{14}$ by
\begin{align}
   \frac1{\log^2y}\sum_{p_1p_2\mid mh}\mu^2(p_1p_2)\log p_1\log p_2d_r\bigg(\frac{mh}{p_1p_2}\bigg)\notag
\end{align}
and
\begin{align}
   \frac1{\log^2y}\sum_{q_1q_2\mid mk}\mu^2(p_1p_2)\log q_1\log q_2d_r\bigg(\frac{mk}{q_1q_2}\bigg)\notag
\end{align}
respectively, the error in calculation of $\mathcal{M}_1(H,T)$ caused by
this is from the terms with $\max(i_1,i_2,j_1,j_2)\ge2$ and the terms
with $p_1=p_2$ or $q_1=q_2$. Since
\begin{align}
   \sum_{i\ge2}\sum_{p^i\le y}\frac{\log^ap}{p^i}=O(1),\notag
\end{align}
for $\forall a\ge0$, we have the sum of the terms with
$\max(i_1,i_2,j_1,j_2)\ge2$
\begin{align}
   \ll&TL^{-3}\sum_{\max(i_1,i_2,j_1,j_2)\ge2}\sum_{m\le y}\frac1m\sum_{h,k\le
   y/m}\frac{1}{hk}\sum_{p_1^{i_1}p_2^{i_2}\mid mh}\log
   p_1\log p_2d_r\bigg(\frac {mh}{p_1^{i_1}p_2^{i_2}}\bigg)\notag\\
   &\cdot\sum_{q_1^{j_1}q_2^{j_2}\mid mk}\log
   q_1\log q_2d_r\bigg(\frac {mk}{q_1^{j_1}q_2^{j_2}}\bigg)\notag\\
   \ll&TL^{-3}\sum_{i_1\ge2}\sum_{p_1^{i_1}\le
   y}\frac{\log p_1}{p_1^i}\sum_{i_2\ge1}\sum_{p_2\le y}\frac{\log p_2}{p_2^{i_2}}\sum_{j_1\ge1}\sum_{q_1\le y}\frac{\log q_1}{q_1^{j_1}}
   \sum_{j_2\ge1}\sum_{q_2\le y}\frac{\log q_2}{q_2^{j_2}}\notag\\
   &\cdot\sum_{m\le y}\frac1m\sum_{h,k\le
   y/m}\frac{d_r(mh)d_r(mk)}{hk}\notag\\
   =&O(TL^{r^2+2r})\notag
\end{align}
and the sum of the terms with $p_1=p_2$
\begin{align}
   \ll&TL^{-3}\sum_{i_1,i_2,j_1,j_2\ge1}\sum_{m\le y}\frac1m\sum_{h,k\le
   y/m}\frac{1}{hk}\sum_{p_1^{i_1+i_2}\mid mh}\log^2
   p_1d_r\bigg(\frac {mh}{p_1^{i_1+i_2}}\bigg)\notag\\
   &\cdot\sum_{q_1^{j_1}q_2^{j_2}\mid mk}\log
   q_1\log q_2d_r\bigg(\frac {mk}{q_1^{j_1}q_2^{j_2}}\bigg)\notag\\
   =&O(TL^{r^2+2r}).\notag
\end{align}
This is also valid to the terms with $q_1=q_2$. So
\begin{align}
   \mathcal {M}_{14}&=T(\log y)^{-4}\sum_{m\le y}\frac{\phi(m)}{m^2}\notag\\
   &\cdot\sum_{h,k\le y/m}\frac{\log\frac{T}{2\pi
   hk}P_2([mh]_y)P_2([mk]_y)}{hk}\sum_{p_1p_2\mid mh}\mu^2(p_1p_2)\log p_1\log p_2\notag\\
   &\cdot d_r\bigg(\frac{mh}{p_1p_2}\bigg)\sum_{q_1q_2\mid mk}\mu^2(q_1q_2)\log q_1\log
   q_2d_r\bigg(\frac{mk}{q_1q_2}\bigg)+O(TL^{r^2+2r}).\notag
\end{align}
We may also replace the sums
\begin{align}
   \sum_{p_1p_2\mid mh}\mu^2(p_1p_2)\log p_1\log p_2d_r\bigg(\frac{mh}{p_1p_2}\bigg),\notag
\end{align}
\begin{align}
   \sum_{q_1q_2\mid mk}\mu^2(q_1q_2)\log q_1\log q_2d_r\bigg(\frac{mk}{q_1q_2}\bigg)\notag
\end{align}
by
\begin{align}
    &\sum_{{i_1'+i_1''=2}\atop{i_1',i_1''\ge0}}C_2^{i_1'}\sum_{p_1\cdots p_{i_1'}\mid m}\mu^2(p_1\cdots p_{i_1'})\log p_1\cdots\log p_{i_1'}\notag\\
    &\times\sum_{p_{i_1'+1}\cdots p_{i_1'+i_1''}\mid h}\log p_{i_1'+1}\cdots\log p_{i_1'+i_1''}d_r\bigg(\frac{mh}{p_1p_2}\bigg),\notag
\end{align}
\begin{align}
   &\sum_{{i_2'+i_2''=2}\atop{i_2',i_2''\ge0}}C_2^{i_2'}\sum_{q_1\cdots q_{i_2'}\mid m}\mu^2(q_1\cdots q_{i_2'})\log q_1\cdots\log q_{i_2'}\notag\\
    &\times\sum_{q_{i_2'+1}\cdots q_{i_2'+i_2''}\mid k}\log q_{i_2'+1}\cdots\log q_{i_2'+i_2''}d_r\bigg(\frac{mk}{q_1q_2}\bigg)\notag
\end{align}
respectively in $\mathcal{M}_{14}$ with an error $O(TL^{r^2+2r})$
as before. Then, we have
\begin{align}
   \mathcal{M}_{14}&=T(\log y)^{-4}\sum_{{i_1'+i_1''=2}\atop{i_1',i_1''\ge0}}\sum_{{i_2'+i_2''=2}\atop{i_2',i_2''\ge0}}C_2^{i_1'}C_2^{i_2'}\sum_{m\le
   y}\frac{\phi(m)}{m^2}\notag\\
   &\sum_{p_1\cdots p_{i_1'}\mid m}\mu^2(p_1\cdots p_{i_1'})\log p_1\cdots\log p_{i_1'}\sum_{q_1\cdots q_{i_2'}\mid m}
   \mu^2(q_1\cdots q_{i_2'})\log q_1\cdots\log q_{i_2'}\notag\\
   &\sum_{h\le
   y/m}\frac{d_r(\frac{mh}{p_1p_2})}{h}P_2([mh]_y)\sum_{p_{i_1'+1}\cdots p_{i_1'+i_1''}\mid h}\log p_{i_1'+1}\cdots\log
   p_{i_1'+i_1''}\notag\\
   &\sum_{k\le y/m}\frac{d_r(\frac{mk}{q_1q_2})}{k}P_2([mk]_y)
   \sum_{q_{i_2'+1}\cdots q_{i_2'+i_2''}\mid k}\log q_{i_2'+1}\cdots\log
   q_{i_2'+i_2''}\log\frac{T}{2\pi hk}\notag
\end{align}
plus an error $O(TL^{r^2+2r})$. We apply Lemma 3.8 to the sum over $h$ and $k$ to obtain
\begin{align}
\label{4.1}
   \mathcal{M}_{14}&=\sum_{{i_1'+i_1''=2}\atop{i_1',i_1''\ge0}}\sum_{{i_2'+i_2''=2}\atop{i_2',i_2''\ge0}}
   \frac{T(\log y)^{2r+i_1''+i_2''-3}C_2^{i_1'}C_2^{i_2'}}{(r+i_1''-1)!(r+i_2''-1)!}
   \sum_{m\le y}\frac{\phi(m)}{m^2}\notag\\
   &\sum_{p_1\cdots p_{i_1'}\mid m}\mu^2(p_1\cdots p_{i_1'})\log p_1\cdots\log
   p_{i_1'}
   \sum_{q_1\cdots q_{i_2'}\mid m}\mu^2(q_1\cdots q_{i_2'})\log q_1\cdots\log
   q_{i_2'}\notag\\
   &\sigma_r\bigg(\frac{m}{p_1\cdots p_{i_1'}}\bigg)\sigma_r\bigg(\frac{m}{q_1\cdots
   q_{i_2'}}\bigg)G([m]_y)+\epsilon_1+\epsilon_2+\epsilon_3+O(TL^{r^2+2r})\tag{4.1}
\end{align}
where
\begin{align}
   G(\alpha)=&\eta^{-1}(1-\alpha)^{2r+i_1''+i_2''}Q_{2,r+i_1''-1}(\alpha)Q_{2,r+i_2''-1}(\alpha)\notag\\
   &-(1-\alpha)^{2r+i_1''+i_2''+1}Q_{2,r+i_1''}(\alpha)Q_{2,r+i_2''-1}(\alpha)\notag\\
   &-(1-\alpha)^{2r+i_1''+i_2''+1}Q_{2,r+i_1''-1}(\alpha)Q_{2,r+i_2''}(\alpha)\notag
\end{align}
and
\begin{align}
   &\epsilon_1\ll TL^{-4}\sum_{m\le
   y}\frac{\sigma_r(m)L^{r+i_1'+i_2'+1}}{m}j_{\tau_0}(m)d_r(m)L^{r+i_1''+i_2''-1}\notag\\
   &\epsilon_2\ll TL^{-4}\sum_{m\le
   y}\frac{j_{\tau_0}(m)d_r(m)L^{r+i_1''+i_2''-1}}{m}\sigma_r(m)L^{r+i_1'+i_2'+1}\notag\\
   &\epsilon_3\ll TL^{-4}\sum_{m\le
   y}\frac{j_{\tau_0}(m)d_r(m)L^{r+i_1''+i_2''}}{m}j_{\tau_0}(m)d_r(m)L^{r+i_1''+i_2''-1}\notag
\end{align}
by (\ref{2.4}), (\ref{2.5}), Lemma 3.1, Lemma 3.2 and an argument as before.
Since
\begin{align}
   |\sigma_r(m)|\ll d_r(m)j_\tau(m) \ \ \ \  for\ \ \ \
   0<\tau\le1\notag
\end{align}
(see (5.13) of N.Ng \cite{Ng}), it follows that
\begin{align}
   \epsilon_1\ll TL^{2r}\sum_{m\le
   y}\frac{d_r(m)^2j_1(m)j_{\tau_0}(m)}{m}\ll TL^{r^2+2r}.\notag
\end{align}
A similar calculation gives $\epsilon_2,\epsilon_3\ll TL^{r^2+2r}$.
Using Lemma 3.2, we have the sum over $m$ in (\ref{4.1}) is
\begin{align}
   &\sum_{m\le y}\frac{\phi(m)}{m^2}\sum_{\tau=0}^{\min(i_1',i_2')}C_{i_1'}^\tau C_{i_2'}^\tau\tau!\notag\\
   &\cdot\sum_{p_1\cdots p_{i_1'+i_2'-\tau}\mid m}\mu^2(p_1\cdots p_{i_1'+i_2'-\tau})\log^2p_1\cdots\log^2p_\tau\log p_{\tau+1}\cdots\log p_{i_1'+i_2'-\tau}
   \notag\\
   &\cdot\sigma_r\bigg(\frac{m}{p_1\cdots p_{i_1'}}\bigg)\sigma_r\bigg(\frac{m}{p_1\cdots
   p_\tau p_{i_1'+1}\cdots p_{i_1'+i_2'-\tau}}\bigg)G([m]_y)\notag\\
   \sim&(\log y)^{i_1'+i_2'}\sum_{\tau=0}^{\min(i_1',i_2')}\frac{r^{i_1'+i_2'-2\tau}a_{r+1}}{(r^2+i_1'+i_2'-1)!}
   \int_0^1\alpha^{r^2+i_1'+i_2'-1}G(\alpha)d\alpha\notag
\end{align}
plus an error $O(L^{-1})$ smaller by Lemma 3.9. Employing this in (\ref{4.1}), we have
\begin{align}
   \mathcal{M}_{14}\sim&\sum_{{i_1'+i_1''=2}\atop{i_1',i_1''\ge0}}\sum_{{i_2'+i_2''=2}\atop{i_2',i_2''\ge0}}\notag\\
   &\frac{T(\log y)^{(r+1)^2}C_2^{i_1'}C_2^{i_2'}a_{r+1}b_r(i_1',i_2')}
   {(r+i_1''-1)!(r+i_2''-1)!(r^2+i_1'+i_2'-1)!}\hat{l}_{2,2}(\eta,r,i_1',i_2',i_1'',i_2'')\notag
\end{align}
with an error $O(TL^{r^2+2r})$. Here, $\hat{l}_{i_1,i_2}(\eta,r,i_1',i_1'',i_2',i_2'')$ is given by (\ref{1.17}). By similar arguments, we can evaluate
$\mathcal{M}_{11}$,
 $\mathcal{M}_{12}$, $\mathcal{M}_{13}$, and have
\begin{align}
\label{4.2}
   \mathcal{M}_1(H,T)\sim&\sum_{i_1=0,2}\sum_{i_2=0,2}\sum_{{i_1'+i_1''=i_1}\atop{i_1',i_1''\ge0}}\sum_{{i_2'+i_2''=i_2}\atop{i_2',i_2''\ge0}}\notag\\
   &\frac{T(\log y)^{(r+1)^2}C_2^{i_1'}C_2^{i_2'}a_{r+1}b_r(i_1',i_2')}
   {(r+i_1''-1)!(r+i_2''-1)!(r^2+i_1'+i_2'-1)!}\hat{l}_{i_1,i_2}(\eta,r,i_1',i_2',i_1'',i_2'').\tag{4.2}
\end{align}
This proves Theorem 1.2.

\section{Evaluation of $I$}
In this section, we will evaluate $I$ in two steps. First, we
apply the lemmas to manipulate $I$ into a suitable form for
evaluation and express $I=I_{1}+I_2+O(yT^{\frac12+\epsilon}+TL^{(r+1)^2})$. Then, we evaluate
$I_{1}$ in section 5.1 and $I_2$ in section 5.2
respectively. Recall that by (\ref{1.23}),
\begin{align}
\label{5.1}
   I=\sum_{k\le
   y}\frac{a(k)}{k}\sum_{j\le\frac{kT}{2\pi}}b(j)e(-j/k)+O(yT^{\frac12+\epsilon}).\tag{5.1}
\end{align}
Using Perron's formula with $c=1+L^{-1}$, the inner sum is
\begin{align}
   \sum_{j\le\frac{kT}{2\pi}}b(j)e(-j/k)=\frac1{2\pi
   i}\int_{c-iT}^{c+iT}\mathcal{Q}^*(s,\alpha,k)\bigg(\frac{kT}{2\pi}\bigg)^s\frac{ds}{s}+O(kT^\epsilon),\notag
\end{align}
where $\mathcal{Q}^*(s,\alpha,k)=\sum_{j=0}^\infty
b(j)j^{-s}e(-j/k)$. Pulling the contour left to
$c_0=\frac12+L^{-1}$, we have
\begin{align}
\label{5.2}
   \sum_{j\le\frac{kT}{2\pi}}b(j)e(-j/k)=&\frac1{2\pi
   i}\bigg(\int_{c-iT}^{c_0-iT}+\int_{c_0-iT}^{c_0+iT}+\int_{c_0+iT}^{c+iT}\bigg)\mathcal{Q}*(s,\alpha,k)\bigg(\frac{kT}{2\pi}\bigg)^s\frac{ds}{s}\notag\\
   &+R_1+R_{1+i\alpha},\tag{5.2}
\end{align}
where $R_u$ is the residue at $s=u$. By Lemma 3.6 the left and
horizontal edges contribute $yT^{\frac12+\epsilon}$. Moreover by
(\ref{1.24}) it follows that
\begin{align}
   \mathcal{Q}^*(s,\alpha,k)=\sum_{h\le
   y}\frac{a(h)}{h^s}\mathcal{Q}(s,\alpha,h/k),\notag
\end{align}
where $\mathcal{Q}(s,\alpha,h/k)$ is defined by (\ref{3.7}). Let
$H=h/(h,k)$, $K=k/(h,k)$, then $\frac hk=\frac HK$ and $(H,K)=1$. We
deduce
\begin{align}
   R_1=\sum_{h\le
   y}a(h)\mathop {\mathrm{res}}\limits_{s=1}\bigg(\mathcal{Q}(s,\alpha,H/K)\bigg(\frac{TK}{2\pi
   H}\bigg)^ss^{-1}\bigg).\notag
\end{align}
By Lemma 3.5(i),
\begin{align}
\label{5.3}
   R_1=&K\sum_{h\le
   y}a(h)\mathop {\mathrm{res}}\limits_{s=1}\bigg(\zeta^2(s)\bigg(\frac{\zeta'}{\zeta}(s-i\alpha)-\mathcal{G}(s,\alpha,K)\bigg)\bigg(\frac
   T{2\pi HK}\bigg)^ss^{-1}\bigg)\notag\\
   =&\frac T{2\pi}\sum_{h\le y}\frac{a(h)}{H}\notag\\
   &\bigg(\big((\zeta'/\zeta)(\bar{\tau})-\mathcal{G}(1,\alpha,K)\big)\log\bigg(\frac{Te^{2\gamma-1}}{2\pi
   HK}\bigg)+\big((\zeta'/\zeta)'(\bar{\tau})-\mathcal{G}'(1,\alpha,K)\big)\bigg),\tag{5.3}
\end{align}
where $\tau=1+i\alpha$ and $\mathcal{G}(s,\alpha,K)$ given by (\ref{3.8}). Similarly, Lemma 3.5(ii) implies
\begin{align}
\label{5.4}
   R_{1+i\alpha}=&\sum_{h\le
   y}a(h)\mathop {\mathrm{res}}\limits_{s=\tau}\bigg(\mathcal{Q}(s,\alpha,H/K)\bigg(\frac{TK}{2\pi
   H}\bigg)^ss^{-1}\bigg)\notag\\
   =&-\frac T{2\pi}\frac{\zeta^2(\tau)}{\tau}\sum_{h\le
   y}\frac{a(h)}{H}\bigg(\frac T{2\pi
   H}\bigg)^{i\alpha}\frac{K\mathcal{R}_K(\tau)}{\phi(K)}.\tag{5.4}
\end{align}
Combining (\ref{5.1}), (\ref{5.2}), (\ref{5.3}) and (\ref{5.4}), we obtain
\begin{align}
   I=&\frac T{2\pi}\sum_{h,k\le
   y}\frac{a(h)a(k)(h,k)}{hk}\bigg(\log\frac{Te^{2\gamma-1}}{2\pi
   HK}\big((\zeta'/\zeta)(\bar{\tau})-\mathcal{G}(1,\alpha,K)\big)+(\zeta'/\zeta)'(\bar{\tau})\notag\\
   &-\mathcal{G}'(1,\alpha,K)-\frac{\zeta^2(\tau)}{\tau}\bigg(\frac
   T{2\pi
   H}\bigg)^{i\alpha}\frac{K\mathcal{R}_K(\tau)}{\phi(K)}\bigg)+O(yT^{\frac12+\epsilon}).\notag
\end{align}
From (\ref{3.11}), we may write
$\mathcal{G}^{(j)}(1,\alpha,K)=\sum_{p\mid
K}p^{i\alpha}\log^{j+1}p+O(C_j(K))$, for $j=0,1$. By Lemma 3.7, the
$O(C_j(K))$ term contributes $O(TL^{(r+1)^2})$. Whence
\begin{align}
   I=&\frac T{2\pi}\sum_{h,k\le
   y}\frac{a(h)a(k)(h,k)}{hk}\bigg(\log\frac{Te^{2\gamma-1}}{2\pi
   HK}\big((\zeta'/\zeta)(\bar{\tau})-\sum_{p\mid
   K}p^{i\alpha}\log p\big)\notag\\
   &+(\zeta'/\zeta)'(\bar{\tau})-\sum_{p\mid
   K}p^{i\alpha}\log^2p-\frac{\zeta^2(\tau)}{\tau}\bigg(\frac T{2\pi
   H}\bigg)^{i\alpha}\frac{K\mathcal{R}_K(\tau)}{\phi(K)}\bigg)\notag
\end{align}
plus an error $O(yT^{\frac12+\epsilon}+TL^{(r+1)^2})$ with $\tau=1+i\alpha$. It follows that
\begin{align}
   I=&\frac T{2\pi}\sum_{h,k\le y}\frac{a(h)a(k)}{hk}\sum_{{m\mid h}\atop{m\mid k}}m\sum_{n\mid
   m}\frac{\mu(n)}{n}\notag\\
   &\cdot\bigg(\log\frac{Te^{2\gamma-1}m^2}{2\pi
   hkn^2}\big((\zeta'/\zeta)(\bar{\tau})-\sum_{p\mid\frac{nk}{m}}p^{i\alpha}\log
   p\big)+(\zeta'/\zeta)'(\bar{\tau})-\sum_{p\mid\frac{nk}{m}}p^{i\alpha}\log^2p\notag\\
   &-\frac{\zeta^2(\tau)}{\tau}\bigg(\frac{Tm}{2\pi
   nh}\bigg)^{i\alpha}\bigg(\frac{nk}{m}\bigg)\frac{\mathcal{R}_{\frac{nk}{m}}(\tau)}{\phi(\frac{nk}{m})}\bigg)+O(yT^{\frac12+\epsilon}+TL^{(r+1)^2}).\notag
\end{align}
by inserting the identity
\begin{align}
   f((h,k))=\sum_{{m\mid h}\atop{m\mid k}}\sum_{n\mid
   m}\mu(n)f\bigg(\frac mn\bigg).\notag
\end{align}
Interchanging summation order and making the variable changes
$h\rightarrow hm$, $k\rightarrow km$ yields
\begin{align}
   I=&\frac T{2\pi}\sum_{m\le y}\frac1m\sum_{n\mid
   m}\frac{\mu(n)}{n}\sum_{h,k\le\frac
   ym}\frac{a(mh)a(mk)}{hk}\notag\\
   &\bigg(\log\frac{Te^{2\gamma-1}}{2\pi
   hkn^2}\big((\zeta'/\zeta)(\bar{\tau})-\sum_{p\mid nk}p^{i\alpha}\log
   p\big)+(\zeta'/\zeta)'(\bar{\tau})\notag\\
   &-\sum_{p\mid
   nk}p^{i\alpha}\log^2p-\frac{\zeta^2(\tau)}{\tau}\bigg(\frac T{2\pi
   nh}\bigg)^{i\alpha}\frac{nk\mathcal{R}_{nk}(\tau)}{\phi(nk)}\bigg)+O(yT^{\frac12+\epsilon}+TL^{(r+1)^2}).\notag
\end{align}
Rearrange this as $I=I_1+I_2+O(yT^{\frac12+\epsilon}+TL^{(r+1)^2})$ where
\begin{align}
   I_1=&\frac T{2\pi}\sum_{m\le y}\frac1m\sum_{n\mid
   m}\frac{\mu(n)}{n}\sum_{h,k\le\frac
   ym}\frac{a(mh)a(mk)}{hk}\notag\\
   &\cdot\bigg(-\log\frac{Te^{2\gamma-1}}{2\pi
   hkn^2}\sum_{p\mid nk}p^{i\alpha}\log
   p-\sum_{p\mid nk}p^{i\alpha}\log^2p\bigg),\notag
\end{align}
\begin{align}
\label{5.5}
   I_2=&\frac T{2\pi}\sum_{m\le y}\frac1m\sum_{n\mid
   m}\frac{\mu(n)}{n}\sum_{h,k\le\frac
   ym}\frac{a(mh)a(mk)}{hk}\notag\\
   &\cdot\bigg(\log\frac{Te^{2\gamma-1}}{2\pi hkn^2}\frac{\zeta'}{\zeta}(\bar{\tau})
   +\bigg(\frac{\zeta'}{\zeta}\bigg)'(\bar{\tau})-\frac{\zeta^2(\tau)}{\tau}\bigg(\frac T{2\pi
   nh}\bigg)^{i\alpha}\frac{nk\mathcal{R}_{nk}(\tau)}{\phi(nk)}\bigg).\tag{5.5}
\end{align}
The first sum is
\begin{align}
   I_1=&\frac T{2\pi}\sum_{m\le y}\frac1m\sum_{n\mid
   m}\frac{\mu(n)}{n}\sum_{h,k\le\frac
   ym}\frac{a(mh)a(mk)}{hk}\notag\\
   &\cdot\bigg(-\log\frac{Te^{2\gamma-1}}{2\pi
   hk}\sum_{p\mid k}p^{i\alpha}\log
   p-\sum_{p\mid k}p^{i\alpha}\log^2p+O(L\log n)\bigg).\notag
\end{align}
A calculation shows that the $O(L\log n)$ contributes
$O(TL^{(r+1)^2})$. Then we deduce that
\begin{align}
\label{5.6}
   I_1=&\frac T{2\pi}\sum_{m\le y}\frac{\phi(m)}{m^2}\sum_{h,k\le\frac
   ym}\frac{a(mh)a(mk)}{hk}\notag\\
   &\cdot\bigg(-\log\frac{Te^{2\gamma-1}}{2\pi
   hk}\sum_{p\mid k}p^{i\alpha}\log
   p-\sum_{p\mid k}p^{i\alpha}\log^2p\bigg)+O(TL^{(r+1)^2}).\tag{5.6}
\end{align}
This puts $I_1$ in a suitable form we need. We now simplify $I_2$ by
substituting the Laurent expansions
\begin{align}
   &(\zeta'/\zeta)(\bar{\tau})=(i\alpha)^{-1}+O(1),\notag\\
   &(\zeta'/\zeta)'(\bar{\tau})=(i\alpha)^{-2}+O(1),\notag\\
   &\zeta^2(\bar{\tau})\tau^{-1}=(i\alpha)^{-2}+(2\gamma-1)(i\alpha)^{-1}+O(1)\notag
\end{align}
in (\ref{5.5}). The $O(1)$ terms of these laurent expansions contribute
\begin{align}
   &TL\sum_{m\le y}\frac1m\sum_{n\mid m}\frac1n\sum_{h,k\le \frac
   ym}\frac{a_r(mh)a_r(mk)}{hk}\ll TL^{(r+1)^2}\notag
\end{align}
by a calculation similar as in section 4, and
\begin{align}
   T\sum_{m\le y}\frac1m\sum_{n\mid m}1\sum_{h\le\frac
   ym}\frac{a(mh)}{h}\bigg|\sum_{k\le\frac ym}a(mk)f(nk)\bigg|\ll
   TL^{r^2+2r}\notag
\end{align}
by Lemma 3.1, Lemma 3.10 and a calculation as before, for $f(k)=\mathcal{R}_k(1+i\alpha)/\phi(k)$
is multiplicative with $f(p^a)\ll p^{-a}$. Thus
we deduce
\begin{align}
\label{5.7}
   I_2=&\frac T{2\pi}\sum_{m\le y}\frac1m\sum_{n\mid
   m}\frac{\mu(n)}{n}\sum_{h,k\le\frac
   ym}\frac{a(mh)a(mk)}{hk}\notag\\
   &\cdot\bigg(\frac{1+i\alpha\log\frac T{2\pi
   hkn^2}-(\frac T{2\pi
   hn})^{i\alpha}\frac{nk\mathcal{R}_{nk}(\tau)}{\phi(nk)}}{(i\alpha)^2}\bigg)+O(TL^{(r+1)^2}).\tag{5.7}
\end{align}
\subsection{Evaluation of $I_1$}
By (\ref{5.6}) it follows that
\begin{align}
   I_1=&\frac T{2\pi}\sum_{m\le y}\frac{\phi(m)}{m^2}\sum_{h,k\le\frac
   ym}\frac{a(mh)a(mk)}{hk}\notag\\
   &\times\bigg(-\log\frac{Te^{2\gamma-1}}{2\pi
   hk}\sum_{p\mid k}p^{i\alpha}\log
   p-\sum_{p\mid k}p^{i\alpha}\log^2p\bigg)+O(TL^{(r+1)^2})\notag
\end{align}
with $a(n)=d_r(n)P_1(\frac{\log n}{\log y})+d_r^*(n)P_2(\frac{\log
n}{\log y})$. As in section 4, we may replace $d_r^*(n)$ by
$\frac1{\log^2y}\sum_{p_1p_2\mid n}\mu^2(p_1p_2)\log p_1\log p_2$,
and the overall error in evaluating $I_1$ caused by this is
$O(L^{-1})$ smaller than the main term, which actually is
$O(TL^{(r+1)^2})$. It follows that
\begin{align}
\label{5.8}
   I_1=\frac
   T{2\pi}\sum_{i_1=0,2}\sum_{i_2=0,2}\big(-L(a_{i_1,i_2,0,0,1})+a_{i_1,i_2,1,0,1}+a_{i_1,i_2,0,1,1}-a_{i_1,i_2,0,0,2}\big)\tag{5.8}
\end{align}
plus an error $O(TL^{(r+1)^2})$, where for $u,v,w\in \mathbb{Z}_{\ge0}$ we define $a_{i_1,i_2,u,v,w}$ to be the
sum
\begin{align}
   &\frac1{(\log y)^{i_1+i_2}}\sum_{m\le
   y}\frac{\phi(m)}{m^2}\sum_{h\le\frac
   ym}\frac{P_{i_1}([mh]_y)\log^uh}{h}\notag\\
   &\cdot\sum_{p_1\cdots p_{i_1}\mid
   mh}\mu^2(p_1\cdots p_{i_1})\log p_1\cdots\log p_{i_1}d_r\bigg(\frac{mh}{p_1\cdots
   p_{i_1}}\bigg)\sum_{k\le\frac
   ym}\frac{P_{i_2}([mk]_y)\log^vk}{k}\notag\\
   &\cdot\sum_{q_1\cdots q_{i_2}\mid
   mk}\mu^2(q_1\cdots q_{i_2})\log q_1\cdots\log q_{i_2}d_r\bigg(\frac{mk}{q_1\cdots
   q_{i_2}}\bigg)\sum_{p\mid k}p^{i\alpha}\log^wp.\notag
\end{align}
Observe that
$a_{i_1,i_2,u,v,w}$
\begin{align}
\label{5.9}
   \sim&\frac1{(\log y)^{i_1+i_2}}\sum_{{i_1'+i_1''=i_1}\atop{i_1',i_1''\ge0}}\sum_{{i_2'+i_2''=i_2}\atop{i_2',i_2''\ge0}}C_{i_1}^{i_1'}C_{i_2}^{i_2'}\sum_{m\le
   y}\frac{\phi(m)}{m^2}\notag\\
   &\cdot\sum_{p_1\cdots p_{i_1'}\mid
   m}\mu^2(p_1\cdots p_{i_1'})\log p_1\cdots\log p_{i_1'}\sum_{q_1\cdots q_{i_2'}\mid
   m}\mu^2(q_1\cdots q_{i_2'})\log q_1\cdots\log q_{i_2'}\notag\\
   &\cdot\sum_{h\le\frac ym}\frac{P_{i_1}([mh]_y)\log^uh}{h}\sum_{p_{i_1'+1}\cdots
   p_{i_1'+i_1''}\mid h}\log p_{i_1'+1}\cdots\log
   p_{i_1'+i_1''}d_r\bigg(\frac{mh}{p_1\cdots p_{i_1}}\bigg)\notag\\
   &\cdot\sum_{k\le\frac
   ym}\frac{P_{i_2}([mk]_y)\log^vk}{k}\sum_{q_{i_2'+1}\cdots
   q_{i_2'+i_2''}\mid k}\log q_{i_2'+1}\cdots\log q_{i_2'+i_2''}d_r\bigg(\frac{mk}{q_1\cdots
   q_{i_2}}\bigg)\notag\\
   &\cdot\sum_{p\mid k}p^{i\alpha}\log^wp\tag{5.9}
\end{align}
plus an error $O(L^{-1})$ smaller. For $p^{i\alpha}=\sum_{j=0}^\infty\frac{(i\alpha)^j}{j!}\log^jp$,
the sum over $k$ in (\ref{5.9}) can be replaced by
\begin{align}
   &\sum_{j=0}^\infty\frac{(i\alpha)^j}{j!}\bigg(\sum_{k\le\frac
   y{m}}\frac{P_{i_2}([mk]_y)(\log k)^v}{k}\notag\\
   &\cdot\sum_{pq_{i_2'+1}\cdots
   q_{i_2'+i_2''}\mid
   k}(\log p)^{j+w}\log q_{i_2'+1}\cdots\log q_{i_2'+i_2''}d_r\bigg(\frac{mk}{pq_1\cdots
   q_{i_2}}\bigg)d_r(p)\notag\\
   &+i_2''\sum_{k\le\frac
   y{m}}\frac{P_{i_2}([mk]_y)(\log k)^v}{k}\notag\\
   &\cdot \sum_{q_{i_2'+1}\cdots
   q_{i_2'+i_2''}\mid
   k}(\log q_{i_2'+1})^{j+w+1}\log q_{i_2'+2}\cdots\log q_{i_2'+i_2''}d_r\bigg(\frac{mk}{q_1\cdots
   q_{i_2}}\bigg)\bigg).\notag
\end{align}
Here, we ignore the terms with $k$ contains square of $p$, for all these terms contribute $O(L^{-1})$ smaller than the main term in the calculation of $a_{i_1,i_2,u,v,w}$.
Substituting this into (\ref{5.9}), we denote
\begin{align}
\label{5.10}
   a_{i_1,i_2,u,v,w}=A_1+A_2\tag{5.10}
\end{align}
plus an error $O(L^{-1})$ smaller with obvious meaning. Then a calculation similar to $\mathcal{M}_1$
in section 4 establishes
\begin{align}
\label{5.11}
   A_1\sim&a_{r+1}\sum_{j=0}^\infty\frac{(i\alpha\log
   y)^j}{j!}(\log y)^{r^2+2r+u+v+w}\notag\\
   &\cdot\sum_{{i_1'+i_1''=i_1}\atop{i_1',i_1''\ge0}}\sum_{{i_2'+i_2''=i_2}\atop{i_2',i_2''\ge0}}rc_r(i_1',i_2',i_1'',i_2'')\beta(w+j,r+i_2'')\notag\\
   &\cdot l_{i_1,i_2}(i_1',i_2',i_1''+u,i_2''+v+w+j)\tag{5.11}
\end{align}
and
\begin{align}
\label{5.12}
   A_2\sim&a_{r+1}\sum_{j=0}^\infty\frac{(i\alpha\log
   y)^j}{j!}(\log y)^{r^2+2r+u+v+w}\notag\\
   &\cdot\sum_{{i_1'+i_1''=i_1}\atop{i_1',i_1''\ge0}}\sum_{{i_2'+i_2''=i_2}\atop{i_2',i_2''\ge0}}i_2''(r+i_2''-1)c_r(i_1',i_2',i_1'',i_2'')\beta(w+j+1,r+i_2''-1)\notag\\
   &\cdot l_{i_1,i_2}(i_1',i_2',i_1''+u,i_2''+v+w+j),\tag{5.12}
\end{align}
with $c_r(i_1',i_2',i_1'',i_2'')$ given by (\ref{1.13}) and $l_{i_1,i_2}(\overrightarrow{n})$ given by (\ref{1.9}). Whence, from
(\ref{5.8}), (\ref{5.10}), (\ref{5.11}) and (\ref{5.12}), we obtain
\begin{align}
\label{5.13}
   I_1\sim&\frac
   T{2\pi}a_{r+1}L^{(r+1)^2+1}\sum_{j=0}^\infty\frac{z^j\eta^{j+(r+1)^2+1}}{j!}\notag\\
   &\cdot\sum_{i_1=0,2}
   \sum_{i_2=0,2}\sum_{{i_1'+i_1''=i_1}\atop{i_1',i_1''\ge0}}\sum_{{i_2'+i_2''=i_2}\atop{i_2',i_2''\ge0}}c_r(i_1',i_2',i_1'',i_2'')\big(rh_{i_1,i_2}(i_1',i_2',i_1''+1,i_2''+1,j)\notag\\
   &+i_2''(r+i_2''-1)h_{i_1,i_2}(i_1',i_2',i_1''+1,i_2'',j+1)\big)\tag{5.13}
\end{align}
with $h_{i_1,i_2}$ given by (\ref{1.11}). Here, we have applied the formula $\beta(a,b)-\beta(a+1,b)=\beta(a,b+1)$ for $\forall a,b\ge1$.
\subsection{Evaluation of $I_2$}
From (\ref{1.6}), we recall that $a(n)=d_r(n)P_1(\frac{\log n}{\log y})+d_r^*(n)P_2(\frac{\log
n}{\log y})$ and use $\frac1{\log^2y}\sum_{p_1p_2\mid
n}\mu^2(p_1p_2)\log p_1\log p_2$ to replace $d_r^*(n)$ in (\ref{5.7}) as
before, then we may denote
\begin{align}
\label{5.14}
   I_2=\frac{T}{2\pi}\sum_{i_1=0,2}\sum_{i_2=0,2}a_{i_1,i_2}'+O(TL^{(r+1)^2})\tag{5.14}
\end{align}
with $a_{i_1,i_2}'$ defined by the sum
\begin{align}
\label{5.15}
   &\frac1{(\log y)^{i_1+i_2}}\sum_{m\le y}\frac1m\sum_{n\mid
   m}\frac{\mu(n)}{n}\notag\\
   &\sum_{h\le\frac
   ym}\frac{P_{i_1}([mh]_y)}{h}\sum_{p_1\cdots p_{i_1}\mid
   mh}\mu^2(p_1\cdots p_{i_1})\log p_1\cdots\log p_{i_1}d_r\bigg(\frac{mh}{p_1\cdots
   p_{i_1}}\bigg)\notag\\
   &\sum_{k\le\frac
   ym}\frac{P_{i_2}([mk]_y)}{k}\sum_{q_1\cdots q_{i_2}\mid
   mk}\mu^2(q_1\cdots q_{i_2})\log q_1\cdots\log q_{i_2}d_r\bigg(\frac{mh}{q_1\cdots
   q_{i_2}}\bigg)\notag\\
   &\bigg(\frac{1+i\alpha\log\frac T{2\pi
   hkn^2}-(\frac T{2\pi
   hn})^{i\alpha}\frac{nk\mathcal{R}_{nk}(\tau)}{\phi(nk)}}{(i\alpha)^2}\bigg)\tag{5.15}
\end{align}
By an argument as before, we have
\begin{align}
\label{5.16}
   a_{i_1,i_2}'&\sim\frac1{(\log y)^{i_1+i_2}}\sum_{{i_1'+i_1''=i_1}\atop{i_1',i_1''\ge0}}\sum_{{i_2'+i_2''=i_2}\atop{i_2',i_2''\ge0}}C_{i_1}^{i_1'}C_{i_2}^{i_2'}\sum_{m\le
   y}\frac{1}{m}\sum_{n\mid m}\frac{\mu(n)}{n}\notag\\
   &\sum_{p_1\cdots p_{i_1'}\mid
   m}\mu^2(p_1\cdots p_{i_1'})\log p_1\cdots\log p_{i_1'}\sum_{q_1\cdots q_{i_2'}\mid
   m}\mu^2(q_1\cdots q_{i_2'})\log q_1\cdots\log q_{i_2'}\notag\\
   &\sum_{h\le\frac ym}\frac{P_{i_1}([mh]_y)}{h}\sum_{p_{i_1'+1}\cdots
   p_{i_1'+i_1''}\mid h}\log p_{i_1'+1}\cdots\log
   p_{i_1'+i_1''}d_r\bigg(\frac{mh}{p_1\cdots p_{i_1}}\bigg)\notag\\
   &\sum_{k\le\frac
   ym}\frac{P_{i_2}([mk]_y)}{k}\sum_{q_{i_2'+1}\cdots
   q_{i_1'+i_2''}\mid k}\log q_{i_1'+1}\cdots\log q_{i_1'+i_2''}d_r\bigg(\frac{mk}{q_1\cdots
   q_{i_2}}\bigg)\notag\\
   &\bigg(\frac{1+i\alpha\log\frac T{2\pi
   hkn^2}-(\frac T{2\pi
   hn})^{i\alpha}\frac{nk\mathcal{R}_{nk}(\tau)}{\phi(nk)}}{(i\alpha)^2}\bigg)\tag{5.16}
\end{align}
plus an error $O(L^{-1})$ smaller. Since all the terms with $k$ that contains square of
$q\in\{q_{i_2'+1},\cdots,q_{i_2'+i_2''}\}$ contribute $O(L^{-1})$ smaller than the main term,
we may ignore these terms in the following argument. Let $k=q_{i_2'+1}\cdots q_{i_2'+i_2''}k'$. For $f(k)=\mathcal{R}_k(1+i\alpha)/\phi(k)$
is multiplicative with $f(p^a)\ll p^{-a}$, we replace
$\frac{\mathcal{R}_{nk}(\tau)}{\phi(nk)}$ by $f(q_{i_2'+1})\cdots
f(q_{i_2'+i_2''})\frac{\mathcal{T}_{nk';r}(\alpha)}{\phi(nk')}$ with an error
\begin{align}
   \ll|\alpha|^{-2}L^r\sum_{m\le y}\frac{d_r(m)}{m}\sum_{n\mid
   m}|\alpha|^{r+1}L^{2r}\frac{d_r(m)j_{\tau_0}(m)}{n^{1-\epsilon}}\ll
   L^{(r+1)^2}\notag
\end{align}
in the calculation of $a_{i_1,i_2}'$ by Lemma 3.1 and Lemma 3.11.
A calculation shows that
$\mathcal{R}_k(1)=\phi(k)/k$ and $\mathcal{R}_k'(1)=-\phi(k)\log k/k$,
thus it follows that
\begin{align}
\label{5.17}
   \frac{\mathcal{T}_{nk';r}(\alpha)}{\phi(nk')}=\frac1{nk}(1-i\alpha\log(nk'))
   +\sum_{j=2}^r\frac{\mathcal{R}_{nk'}^{(j)}(1)(i\alpha)^j}{\phi(nk')j!}\tag{5.17}
\end{align}
and
\begin{align}
\label{5.18}
   f(p)=&\frac{\mathcal{R}_p(\tau)}{\phi(p)}=\frac1{\phi(p)}(2-p^{i\alpha}-\frac1{p^{1+i\alpha}})\notag\\
   =&\frac1p(1-i\alpha\log
   p)-\frac1p\sum_{j=2}^\infty\frac{(i\alpha)^j\log^jp}{j!}+O\bigg(\frac{\log^2p}{p^2}(i\alpha)^2\bigg)\notag\\
   =&\frac1p\sum_{j=0}^\infty\frac{\Delta(j)(i\alpha)^j\log^jp}{j!}+O\bigg(\frac{\log^2p}{p^2}(i\alpha)^2\bigg)\tag{5.18}
\end{align}
with $\Delta(j)$ given by (\ref{1.15}). Here, the $O(\frac{\log^2p}{p^2}(i\alpha)^2)$ contributes
$O(L^{(r+1)^2})$ in $a_{i_1,i_2}'$ by a calculation as before. Substituting (\ref{5.17}) and (\ref{5.18}) into
(\ref{5.16}), we have the expressing within the brackets of (\ref{5.16})
simplifies to
\begin{align}
   -\sum_{u+j+j_1+\cdots+j_{i_2''}\ge2}\frac{(\log\frac{T}{2\pi
   hn})^u\mathcal{R}_{nk'}^{(j)}(1)(i\alpha)^{u+j}\prod_{\tau=1}^{i_2''}\bigg((i\alpha)^{j_\tau}\Delta(j_\tau)\log^{j_\tau}q_{i_2'+\tau}\bigg)}{\phi(nk')u!j!j_1!\cdots
   j_{i_2''}!}\notag
\end{align}
by replacing $(\frac T{2\pi hn})^{i\alpha}$ with $\sum_{u=0}^\infty\frac{(i\alpha)^u}{u!}(\log\frac T{2\pi hn})^u$. Employing this in (\ref{5.16}), we have $a_{i_1,i_2}'$ equal to
\begin{align}
   \label{5.181}
   &\frac{-1}{(\log
   y)^{i_1+i_2}}\sum_{{i_1'+i_1''=2}\atop{i_1',i_1''\ge0}}\sum_{{i_2'+i_2''=2}\atop{i_2',i_2''\ge0}}C_{i_1}^{i_1'}C_{i_2}^{i_2'}\notag\\
   &\cdot\sum_{u+j+j_1+\cdots+j_{i_2''}\ge2}\frac{(i\alpha)^{u+j+j_1+\cdots+j_{i_2''}}\Delta(j_1)\cdots\Delta(j_{i_2''})}{u!j!j_1!\cdots
   j_{i_2''}!}\sum_{m\le
   y}\frac{1}{m}\sum_{n\mid m}\frac{\mu(n)}{n}\notag\\
   &\cdot\sum_{p_1\cdots p_{i_1'}\mid
   m}\mu^2(p_1\cdots p_{i_1'})\log p_1\cdots\log p_{i_1'}\sum_{q_1\cdots q_{i_2'}\mid
   m}\mu^2(q_1\cdots q_{i_2'})\log q_1\cdots\log q_{i_2'}\notag\\
   &\cdot\sum_{h\le\frac ym}\frac{P_{i_1}([mh]_y)}{h}\sum_{p_{i_1'+1}\cdots
   p_{i_1'+i_1''}\mid h}\log p_{i_1'+1}\cdots\log
   p_{i_1'+i_1''}d_r\bigg(\frac{mh}{p_1\cdots p_{i_1}}\bigg)\notag\\
   &\cdot \bigg(\log\frac{T}{2\pi hn}\bigg)^u\sum_{k\le\frac
   ym}\frac{P_{i_2}([mk]_y)}{k}\sum_{q_{i_2'+1}\cdots
   q_{i_2'+i_2''}\mid k}\log^{j_1+1}q_{i_2'+1}\cdots\log^{j_{i_2''+1}}q_{i_2'+i_2''}\notag\\
   &\cdot d_r\bigg(\frac{mk}{q_1\cdots
   q_{i_2}}\bigg)\frac{nk'\mathcal{R}_{nk'}^{(j)}(1)}{\phi(nk')}\tag{5.19}
\end{align}
plus an error $O(L^{(r+1)^2})$. Whence, by Lemma 3.8 and Lemma 3.11, we have the sum over $h$ in (\ref{5.181}) equal to
\begin{align}
   \frac{\sigma_r(\frac{m}{p_1\cdots p_{i_1'}})(\log y)^{r+i_1''+u}}{(r+i_1''-1)!}\int_0^{1-[m]_y}\theta_1^{r+i_1''-1}(\eta^{-1}-\theta_1)^uP_{i_1}([m]_y+\theta_1)d\theta_1\notag
\end{align}
and the sum over $k$ equal to
\begin{align}
   \frac{\sigma_r(\frac{m}{q_1\cdots q_{i_2'}})(-1)^jC_r^jj_1!\cdots j_{i_2''}!(\log\frac ym)^{r+i_2''+j+j_1+\cdots+ j_{i_2''}}}{(r+i_2''+j+j_1+\cdots+ j_{i_2''}-1)!}Q_{i_2,r+i_2''+j+j_1+\cdots+ j_{i_2''}-1}([m]_y)\notag
\end{align}
Employing these into (\ref{5.181}), we interchange the order of the sum and the integration, and by Lemma 3.9, we have that
\begin{align}
\label{5.19}
   a_{i_1,i_2}'\sim&(\log y)^{(r+1)^2}\sum_{{i_1'+i_1''=i_1}\atop{i_1',i_1''\ge0}}\sum_{{i_2'+i_2''=i_2}\atop{i_2',i_2''\ge0}}
   \sum_{u+j+j_1+\cdots+j_{i_2''}\ge2}(i\alpha\log
   y)^{u+j+j_1+\cdots+j_{i_2''}-2}\notag\\
   &\cdot\frac{C_{i_1}^{i_1'}C_{i_2}^{i_2'}(-1)^{j+1}C_r^j\Delta(j_1)\cdots\Delta(j_{i_2''})a_{r+1}b_r(i_1',i_2')}
   {u!(r+i_1''-1)!(r+i_2''+j+j_1+\cdots+j_{i_2''}-1)!(r^2+i_1'+i_2'-1)!}\notag\\
   &\cdot k_{i_1,i_2}(i_1'',u,i_1'+i_2',i_2''+j+j_1+\cdots+j_{i_2''})\tag{5.20}
\end{align}
plus an error $O(L^{(r+1)^2})$. We now simplify the expression of $a_{i_1,i_2}'$ in
three cases.

Case 1. $i_2''=0$. Replacing $j-2$ by $n$ and $j+u-2$ by $j$ in
(\ref{5.19}) respectively, we have
\begin{align}
\label{5.10}
   a_{i_1,i_2}'\sim&(\log y)^{(r+1)^2}\sum_{{i_1'+i_1''=i_1}\atop{i_1',i_1''\ge0}}\sum_{{i_2'+i_2''=i_2}\atop{i_2',i_2''\ge0}}C_{i_1}^{i_1'}C_{i_2}^{i_2'}
   \sum_{j=0}^\infty(i\alpha\log
   y)^{j}\notag\\
   &\cdot\sum_{n=-2}^{\min(r-2,j)}\frac{a_{r+1}b_r(i_1',i_2')(-1)^{n+1}C_r^{n+2}}
   {(j-n)!(r+i_1''-1)!(r+i_2''+n+1)!(r^2+i_1'+i_2'-1)!}\notag\\
   &\cdot k_{i_1,i_2}(i_1'',j-n,i_1'+i_2',i_2''+n+2)\tag{5.21}
\end{align}

Case 2. $i_2''=1$. We have
\begin{align}
\label{5.21}
   a_{i_1,i_2}'\sim&(\log y)^{(r+1)^2}\sum_{{i_1'+i_1''=i_1}\atop{i_1',i_1''\ge0}}\sum_{{i_2'+i_2''=i_2}\atop{i_2',i_2''\ge0}}C_{i_1}^{i_1'}C_{i_2}^{i_2'}
   \sum_{j=0}^\infty(i\alpha\log
   y)^{j}\notag\\
   &\cdot\sum_{n=-2}^{j}\frac{a_{r+1}b_r(i_1',i_2')\sum_{j'=-2}^{\min(r-2,n)}(-1)^{j'+1}C_r^{j'+2}\Delta(n-j')}
   {(j-n)!(r+i_1''-1)!(r+i_2''+n+1)!(r^2+i_1'+i_2'-1)!}\notag\\
   &\cdot k_{i_1,i_2}(i_1'',j-n,i_1'+i_2',i_2''+n+2)\tag{5.22}
\end{align}
by replacing $j-2$ with $j'$, $j+j_1-2$ with $n$ and $j+u+j_1-2$
with $j$ in (\ref{5.19}) respectively.

Case 3. $i_2''=2$. We have
\begin{align}
\label{5.22}
   a_{i_1,i_2}'\sim&(\log y)^{(r+1)^2}\sum_{{i_1'+i_1''=i_1}\atop{i_1',i_1''\ge0}}\sum_{{i_2'+i_2''=i_2}\atop{i_2',i_2''\ge0}}C_{i_1}^{i_1'}C_{i_2}^{i_2'}
   \sum_{j=0}^\infty(i\alpha\log
   y)^{j}\notag\\
   &\cdot\sum_{n=-2}^{j}\frac{a_{r+1}b_r(i_1',i_2')\sum_{j'=-2}^{\min(r-2,n)}(-1)^{j'+1}C_r^{j'+2}\sum_{j_1+j_2=n-j'}\Delta(j_1)\Delta(j_2)}
   {(j-n)!(r+i_1''-1)!(r+i_2''+n+1)!(r^2+i_1'+i_2'-1)!}\notag\\
   &\cdot k_{i_1,i_2}(i_1'',j-n,i_1'+i_2',i_2''+n+2)\tag{5.23}
\end{align}
by replacing $j-2$ with $j'$, $j+j_1+j_2-2$ with $n$ and
$j+u+j_1+j_2-2$ with $j$ in (\ref{5.19}) respectively. Since
\begin{align}
   \sum_{j=0}^{r}(-1)^{j+1}C_r^{j}P(j)=0\notag
\end{align}
for any polynomial $P(j)$ on $j$, we have
\begin{align}
   \sum_{j'=-2}^{\min(r-2,n)}(-1)^{j'+1}C_r^{j'+2}\Delta(n-j')=0\ \ \
   \ \ \ for\ \ \ \ \ \ n>r-2,\notag
\end{align}
and
\begin{align}
   \sum_{j'=-2}^{\min(r-2,n)}(-1)^{j'+1}C_r^{j'+2}\sum_{j_1+j_2=n-j'}\Delta(j_1)\Delta(j_2)=0\ \ \
   \ for\ \ \ \ n>r-2.\notag
\end{align}
So, we simplify the expression of $a_{i_1,i_2}'$ for all case to
\begin{align}
\label{5.23}
   a_{i_1,i_2}'\sim&(\log y)^{(r+1)^2}a_{r+1}\sum_{j=0}^\infty(i\alpha\log
   y)^{j}\sum_{{i_1'+i_1''=i_1}\atop{i_1',i_1''\ge0}}\sum_{{i_2'+i_2''=i_2}\atop{i_2',i_2''\ge0}}c_r(i_1',i_2',i_1'',i_2'')\notag\\
   &\cdot\sum_{n=-2}^{\min(r-2,j)}\frac{\Omega_{r}(i_2'',n)(r+i_2''-1)!}
   {(j-n)!(r+i_2''+n+1)!}
    k_{i_1,i_2}(i_1'',j-n,i_1'+i_2',i_2''+n+2).\tag{5.24}
\end{align}
with $\Omega_{r}(i_2'',n)$ given by (\ref{1.14}). Thus, substituting (\ref{5.23}) into (\ref{5.14}), we have
\begin{align}
\label{5.24}
   I_2=&\frac{T}{2\pi}L^{(r+1)^2+1}a_{r+1}\sum_{j=0}^\infty(z)^{j}\eta^{j+(r+1)^2+1}\sum_{i_1=0,2}\sum_{i_2=0,2}
   \sum_{{i_1'+i_1''=i_1}\atop{i_1',i_1''\ge0}}\sum_{{i_2'+i_2''=i_2}\atop{i_2',i_2''\ge0}}c_r(i_1',i_2',i_1'',i_2'')\notag\\
   &\cdot\sum_{n=-2}^{\min(r-2,j)}\frac{\Omega_{r}(i_2'',n)(r+i_2''-1)!}
   {(j-n)!(r+i_2''+n+1)!}
    k_{i_1,i_2}(i_1'',j-n,i_1'+i_2',i_2''+n+2)\tag{5.25}
\end{align}
plus an error $O(TL^{(r+1)^2})$. Theorem 1.3 follows from (\ref{1.22}), (\ref{4.2}),(\ref{5.13}) and (\ref{5.24}).

\end{document}